# POSTERIOR CONVERGENCE RATES OF DIRICHLET MIXTURES AT SMOOTH DENSITIES

## By Subhashis Ghosal[1] and Aad van der Vaart

### *North Carolina State University and Vrije Universiteit, Amsterdam*

We study the rates of convergence of the posterior distribution for Bayesian density estimation with Dirichlet mixtures of normal distributions as the prior. The true density is assumed to be twice continuously differentiable. The bandwidth is given a sequence of priors which is obtained by scaling a single prior by an appropriate order. In order to handle this problem, we derive a new general rate theorem by considering a countable covering of the parameter space whose prior probabilities satisfy a summability condition together with certain individual bounds on the Hellinger metric entropy. We apply this new general theorem on posterior convergence rates by computing bounds for Hellinger (bracketing) entropy numbers for the involved class of densities, the error in the approximation of a smooth density by normal mixtures and the concentration rate of the prior. The best obtainable rate of convergence of the posterior turns out to be equivalent to the well-known frequentist rate for integrated mean squared error $n^{-2/5}$ up to a logarithmic factor.

**1. Introduction.** Kernel methods for density estimation have been in use for nearly fifty years. Bayesian kernel density estimation using a Dirichlet process on the mixing distribution has been considered more recently (cf. [5, 12]), where the density is viewed as a mixture of normals with an arbitrary mixing distribution and a Dirichlet process (cf. [4]) is used as a prior on the mixing distribution. Efficient Gibbs sampling algorithms for the computation of the posterior based on a Dirichlet mixture process have been developed; see, for instance, [3]. Under certain conditions, posterior consistency of such a Dirichlet mixture prior with a normal kernel has been obtained

Received August 2005; revised April 2006.

[1]Supported in part by NSF Grant DMS-03-49111 and funding from Vrije Universiteit for a visit.

by Ghosal, Ghosh and Ramamoorthi [7]. Ghosal and van der Vaart [9] obtained rates of convergence of the posterior for the Dirichlet mixture in the case that the true density is a location or location-scale mixture of normals with standard deviations bounded away from zero and infinity. Under natural conditions on the prior, they showed that the posterior converges at rate $(\log n)^\kappa / \sqrt{n}$, where $\kappa$ depends on the tail behavior of the base measure of the Dirichlet process. The rate of convergence was obtained by finding a sharp entropy estimate and prior concentration rate for this problem and then applying the general posterior convergence rate theorem of Ghosal, Ghosh and van der Vaart [8]. The fast rate of convergence $(\log n)^\kappa / \sqrt{n}$ arises because a mixture of normals with standard deviations bounded by two positive numbers is "super-smooth." Super-smooth densities can be approximated by kernel estimators with a bandwidth that approaches zero at a logarithmic rate and super-smooth mixtures can be well approximated by finite normal mixtures with a small number of components (cf. Lemma 3.1 of [9]). This leads to small entropy numbers and high prior concentration (comparable to those of finite-dimensional models) with a nearly parametric convergence rate as a consequence. As a consequence of entropy bounds for normal mixtures, Ghosal and van der Vaart [9] also obtained essentially the same convergence rate $(\log n) / \sqrt{n}$ for sieved maximum likelihood estimators (MLE). Under the same super-smoothness condition, Genovese and Wasserman [6] earlier obtained the much weaker convergence rate $n^{-1/6} (\log n)^{(1+\delta)/6}$ for some $\delta > 0$ for sieved MLEs based on Gaussian mixtures.

While it is interesting to observe nearly parametric rates of convergence, the super-smoothness of the true density with a bounded known range for the standard deviation is a restrictive assumption. Scricciolo [13] considered the situation where the true density is still super-smooth, but the prior for the bandwidth parameter contains zero in its support. The resulting rate of convergence is much slower in this case and depends on the decay rate of the prior for the bandwidth at zero. In this paper, we consider the more realistic situation where the density of the observations is smooth, but may not be a mixture of normal densities. A smooth density can be approximated by mixtures of normals, but it is necessary to let the bandwidth (standard deviations of the components) tend to zero and allow an increasing number of components. This increases the complexity of the model and leads to larger entropy and smaller prior concentration, with, as a consequence, a slower rate of convergence of the posterior distribution.

More specifically, we assume that the density of the observations is twice continuously differentiable. Under some regularity conditions, the optimal rate of convergence of a kernel estimator is then $n^{-2/5}$. The main purpose of this paper is to establish that the posterior distribution based on a Dirichlet mixture of normal prior attains the same rate, up to a logarithmic factor. In addition, we obtain the same rate for the sieved maximum likelihood estimator using a sieve consisting of normal mixtures. It may be noted that, even though the estimation of a smooth density is a considerably harder problem than that of a super-smooth density, our obtained rate, which is nearly optimal for the given problem, is still much better than the $n^{-1/6}$ rate Genovese and Wasserman [6] obtained for sieved MLEs in the super-smooth case.

1.1. *Notation.* Throughout the paper $X_1, X_2, \ldots$ are independent and identically distributed (i.i.d.) as $p_0$ on $\mathbb{R}$. The corresponding probability measure is denoted by $P_0$.

The supremum and $L_1$-norm are denoted by $\|\cdot\|_\infty$ and $\|\cdot\|_1$, respectively. For two density functions $f, g : \mathbb{R} \to [0, \infty)$, we let $h$ denote the Hellinger distance defined by $h^2(f, g) = \int (f^{1/2} - g^{1/2})^2 \, d\lambda$, where $\lambda$ is the Lebesgue measure on $\mathbb{R}$. The $\varepsilon$-covering number $N(\varepsilon, S, d)$ of a semi-metric space $S$ relative to the semi-metric $d$ is the minimal number of balls of radius $\varepsilon$ needed to cover $S$. Similarly, the $\varepsilon$-bracketing number $N_{[\cdot]}(\varepsilon, S, d)$ is the minimal number of $\varepsilon$-brackets $[f, g] = \{u : f \le u \le g\}$ needed to cover $S$, the size of a bracket $[f, g]$ being the distance $d(f, g)$ between upper and lower brackets (cf., e.g., [14]). The logarithms of the covering and bracketing numbers are referred to as entropies without and with bracketing.

We write "$\lesssim$" for inequality up to a constant multiple, where the constant is universal or (at least) unimportant for our purposes. An expression $x^{a+}$ in a statement means that the statement holds for $x^{a'}$ for any $a' > a$. Let $\phi(x) = (2\pi)^{-1/2} \exp(-x^2/2)$, the standard normal density, and let $\phi_\sigma(x) = \sigma^{-1}\phi(x/\sigma)$. An asterisk denotes convolution and $p_{F,\sigma} = F * \phi_\sigma$ is a Gaussian mixture with mixing distribution $F$. The distribution which is degenerate at $\theta$ is denoted by $\delta_\theta$. The support of a density $p$ is denoted by $\mathrm{supp}(p)$.

1.2. *Assumptions.* Throughout the paper, we assume that $h(p_0, p_0 * \phi_\sigma) = O(\sigma^2)$ as $\sigma \to 0$. If $p_0$ is twice continuously differentiable with $\int (p_0''/p_0)^2 \times p_0 \, d\lambda < \infty$ and $\int (p_0'/p_0)^4 p_0 \, d\lambda < \infty$, then the condition holds (cf. Lemma 4).

1.3. *Organization.* The main results of the paper are on the convergence rate of the posterior distribution and these are presented in Section 2. The proofs of the main theorems are contained in Sections 9 and 10, and are based on estimates of the entropies of normal mixtures obtained in Section 5, approximation lemmas given in Section 6 and lower bounds on Dirichlet probabilities obtained in Section 7. A general result on posterior convergence rates is obtained in Section 4, which is subsequently used in the proof of the main result in Section 2. The entropy estimates also have applications to rates of convergence of sieved MLEs and posterior distributions relative to sieved priors, as noted in Section 3. The proofs of the theorems in Section 3 are given in Section 11. Sections 4–8 may be of some independent interest.

**2. Main results.** We consider the sequence of priors $\Pi_n$ for $p$ defined structurally as follows:

- $p_{F,\sigma}(x) = \int \phi_\sigma(x - z)\, dF(z)$;
- $F \sim D_\alpha$, the Dirichlet process with base measure $\alpha = \alpha(\mathbb{R})\bar{\alpha}$, where $0 < \alpha(\mathbb{R}) < \infty$ and $\bar{\alpha}$ is a probability measure;
- $\sigma/\sigma_n \sim G$, where $\sigma_n$ is a sequence of positive real numbers converging to zero with $n^{-a_1} \lesssim \sigma_n \lesssim n^{-a_2}$ for some $0 < a_2 < a_1 < 1$ and $G$ is a fixed probability distribution on $(0, \infty)$ satisfying $G(s) \lesssim e^{-\beta s^{-\gamma}}$ as $s \to 0$ and $1 - G(s) \lesssim e^{-\beta s^\gamma}$ as $s \to \infty$ for some $\gamma > 1$ and $\beta > 0$.
- $F$ and $\sigma$ are independent.

If $\Pi_n(p: d(p, p_0) > M\varepsilon_n | X_1, \ldots, X_n) \to 0$ in $P_0^n$-probability for some $M > 0$, we say that $\varepsilon_n \to 0$ is (an upper bound for) the posterior convergence rate relative to $d$.

The proof of the following posterior convergence theorem is given in Section 9.

THEOREM 1. *Suppose that $p_0$ has compact support and that $a_2 \geq (4 + \gamma)^{-1}$. If the base measure $\alpha$ has a continuous and positive density on an interval containing* $\mathrm{supp}(p_0)$, *then the posterior rate of convergence relative to $h$ is*

$$(2.1) \qquad \varepsilon_n = \max\{(n\sigma_n)^{-1/2}(\log n), n^{-1/2}(\sigma_n^{-1})^{(\gamma/2(\gamma-1))+}, \sigma_n^2 \log n\}.$$

*If* $\mathrm{supp}(p_0)$ *is a finite union of intervals and every interval $I$ in the support satisfies $P_0(I) \gtrsim \lambda(I)^a$ for some $a > 0$, then this can be improved to the rate*

$$(2.2) \qquad \varepsilon_n = \max\{(n\sigma_n)^{-1/2}(\log n), n^{-1/2}(\sigma_n^{-1})^{(\gamma/2(\gamma-1))+}, \sigma_n^2\}.$$

*Further, when $G$ is compactly supported, the middle terms on the right-hand side of* (2.1) *and* (2.2) *can be omitted.*

The best rate $n^{-2/5}(\log n)^{4/5}$ in the preceding theorem is obtained in the second assertion with $\gamma = \infty$ (i.e., $G$ is compactly supported) if $\sigma_n$ is chosen to be $n^{-1/5}(\log n)^{2/5}$, nearly equal to the optimal frequentist bandwidth choice $n^{-1/5}$.

A common practice is to consider an inverse gamma prior on $\sigma^2$, which leads to conditional conjugacy and hence to an efficient Gibbs sampling procedure. Unfortunately, our theorem does not apply to this prior, because the inverse gamma prior has a polynomially decaying tail near infinity. Indeed, even with faster-than-exponential decay, the theorem indicates that rates may suffer whenever the support of the prior is noncompact. Because these rates are only upper bounds, a negative conclusion cannot be reached based on these. However, it may be mentioned that even the issue of consistency is open for the inverse gamma prior unless an upper truncation is used (cf. [7]). On the other hand, for a truncated inverse gamma prior, a nearly optimal convergence rate is obtained from Theorem 1, while Gibbs sampling can be implemented easily with an additional acceptance–rejection step to take care of the truncation.

The preceding theorem will be obtained by applying the general posterior convergence rate theorem in Section 4. Estimates of entropy and prior concentration rate obtained in [9] for the super-smooth case will be refined in a way suitable to the present situation in Section 5.

The assumption that $p_0$ is compactly supported is restrictive, in particular in combination with the assumption that $h(p_0, p_0 * \phi_\sigma) = O(\sigma^2)$, which forces $p_0$ to tend to zero smoothly at the boundary points of its support. We do not know if the assumption of compact support can be completely removed, but we note the following extensions of the preceding theorem, which increase the applicability considerably.

Given a smooth function $w : \mathbb{R} \to [0, 1]$ with compact support, we can form a reduced data set $\bar{X}_1, \ldots, \bar{X}_{\bar{n}}$ by rejecting each $X_i$ independently with probability $1 - w(X_i)$, giving a sample from the density $\bar{p}_0 = p_0 w / \int p_0 w \, d\lambda$. The size $\bar{n}$ is distributed binomially with parameters $n$ and $\int p_0 w \, d\lambda$, whence $\bar{n}/n \to \int p_0 w \, d\lambda$ a.s. Conditionally on $\bar{n}$, Theorem 2.1 of [8] can be applied to conclude that the posterior concentrates on Hellinger balls of radius $\varepsilon_{\bar{n}}$ around $\bar{p}_0$. If we choose $w$ to be equal to 1 on a given compact then $\bar{p}_0$ and $p_0$ are proportional on this compact and hence the posterior essentially gives the (conditional) density of the original observations on this compact.

This construction may be appropriate for Bayesian estimation of heavy-tailed densities, but it does change the posterior distribution. Even though we may expect that the change on an interval where $w$ is identically one is minimal, it appears to be difficult to bound the difference. This difficulty can be avoided by applying the preceding with a sequence of truncation functions. For densities with exponentially decreasing tails, it yields a rate of convergence of the posterior relative to the Hellinger (semi-)distance on compact intervals given by $h_k^2(p, q) = \int_{-k}^{k} (p^{1/2} - q^{1/2})^2 \, d\lambda$. For simplicity, in this result we assume that $G$ is compactly supported in $(0, \infty)$. The proof of the following theorem is contained in Section 10.

THEOREM 2. *Suppose that $p_0$ satisfies $P_0[-a, a]^c \le e^{-ca^\gamma}$ for some positive numbers $c$ and $\gamma$, and is twice continuously differentiable with $\int (p_0''/p_0)^2 \times p_0 \, d\lambda < \infty$ and $\int (p_0'/p_0)^4 p_0 \, d\lambda < \infty$. If the base measure $\alpha$ has a continuous and positive density $\alpha'$ satisfying $\alpha'(t) \gtrsim e^{-dt^\gamma}$ for sufficiently large $|t|$, for some positive constant $d$, then the rate of convergence relative to the semi-distance $h_k$ is at least*

$$(2.3) \qquad \varepsilon_n = \max\{(n\sigma_n)^{-1/2}(\log n)^{1+\gamma/2}, \sigma_n^2 \log n\}.$$

**3. Sieve maximum likelihood and sieve priors.** As a byproduct of the upper bounds on the entropy of the set of normal mixtures (necessary for the proofs of our main results), we can also obtain the rate of convergence of sieved MLEs for normal mixtures. We consider sieves of the types

$$(3.1) \qquad \mathcal{P}_n = \{p_{F,\sigma} : F[-a_n, a_n] = 1, b_1\sigma_n \le \sigma \le b_2\sigma_n\},$$

$$(3.2) \qquad \mathcal{P}_n = \{p_{F,\sigma} : F[-a, a]^c \le A(a) \text{ for all } a > 0, b_1\sigma_n \le \sigma \le b_2\sigma_n\}.$$

Here, $a_n$ and $\sigma_n$ are positive sequences and $A : (0, \infty) \to [0, 1]$ is decreasing. We define the sieved MLE as $\hat{p}_n = \arg\max\{\prod_{i=1}^n p(X_i) : p \in \mathcal{P}_n\}$.

The rate of convergence of sieved MLEs relative to $h$ can be obtained from Theorem 4 of [16] or Theorem 3.4.4 of [14]. There is a trade-off between the complexity of the model $\mathcal{P}_n$ and the distance of $\mathcal{P}_n$ to $p_0$. Under the assumption of Section 1.2, the approximation rate is $O(\sigma_n^2)$. The complexity of the model $\mathcal{P}_n$ can be measured through its bracketing entropy. The rate of convergence is the maximum of the approximation error and the solution $\varepsilon_n$ to the equation

$$(3.3) \qquad \int_0^{\varepsilon_n} \sqrt{\log N_{[\cdot]}(\varepsilon, \mathcal{P}_n, h)} \, d\varepsilon \sim \sqrt{n}\varepsilon_n^2.$$

THEOREM 3. *Let $\sigma_n \to 0$ and $a_n \ge e$ so that $\log n \lesssim \log(a_n/\sigma_n) \lesssim \log n$, and let $\hat{p}_n$ be the sieved MLE relative to $\mathcal{P}_n$ given by (3.1). If $p_0$ has compact support and $[-a_n, a_n] \supset \mathrm{supp}(p_0)$ for all sufficiently large $n$, then $h(\hat{p}_n, p_0) = O_P(\varepsilon_n)$ for*

$$(3.4) \qquad \varepsilon_n = \max\{(n\sigma_n)^{-1/2}a_n\log n, \sigma_n^2\}.$$

The apparently best rate $n^{-2/5}(\log n)^{4/5}$ is obtained when $a_n$ is bounded, but $[-a_n, a_n] \supset \mathrm{supp}(p_0)$ and $\sigma_n \sim n^{-1/5}(\log n)^{2/5}$. The optimal order of bandwidth for the classical kernel estimator $\sigma_n \sim n^{-1/5}$ leads to a slightly larger error rate $n^{-2/5}\log n$. Admittedly, these are only upper bounds. In particular, the logarithmic factor may not be sharp.

When $\mathrm{supp}(p_0)$ is not compact, but $p_0/(p_0 * \phi_{\sigma_n})$ are uniformly bounded, we can use the sieves (3.2) to derive the convergence rate. The condition holds, for instance, if $p_0$ is increasing on $(-\infty, a]$, bounded below on $[a, b]$ and decreasing on $[b, \infty)$ for some $a < b$ (cf. Lemma 6 in Section 6).

THEOREM 4. *Let $\sigma_n \to 0$ so that $\log n \lesssim \log \sigma_n^{-1} \lesssim \log n$ and let $\hat{p}_n$ be the sieved MLE relative to $\mathcal{P}_n$ given by (3.2) with $A(a) = e^{-da^{1/\delta}}$, $d, \delta > 0$ constants. If $P_0[-a, a]^c \le A(a)$ for every $a > 0$, then $h(\hat{p}_n, p_0) = O_p(\varepsilon_n)$ for*

$$(3.5) \qquad \varepsilon_n = \max\{(n\sigma_n)^{-1/2}(\log n)^{1+(1\vee 2\delta)/4}, \sigma_n^2\}.$$

The proofs of the theorems in this section are given in Section 11.

The sieves $\mathcal{P}_n$ in (3.1) and (3.2) of the preceding section can also be used to construct a prior for which the posterior converges at the same rate $\varepsilon_n$ as obtained in Theorems 3 and 4. As in Theorem 3.1 of [8], take a minimal collection of Hellinger $\varepsilon_n$-brackets that cover $\mathcal{P}_n$. Consider the uniform prior $\Pi_n$ on the renormalized upper brackets. Then the resulting posterior converges at the rate $\varepsilon_n$.

## 4. A general result on posterior convergence rates.

When the prior $G$ on $\sigma/\sigma_n$ is not compactly supported, existing results on posterior convergence rates (such as Theorem 2.1 of [8]) do not seem to suffice in deriving the rate. Below we obtain a posterior convergence rate theorem where we use a countable decomposition of the space of densities together with conditions on their prior probabilities and entropy numbers with respect to $h$.

Let $X_1, X_2, \ldots$ be i.i.d. with density $p \in \mathcal{P}$. Let $\Pi_n$ be a sequence of priors on $\mathcal{P}$ and let $p_0$ and $P_0$ stand for the true density and the true probability measure, respectively. Let $d$ be a metric which induces convex balls and is bounded above on $\mathcal{P}$ by a multiple of $h$.

THEOREM 5. *Suppose that $\mathcal{P}_n \subset \mathcal{P}$ is such that $\Pi_n(\mathcal{P}_n^c | X_1, \ldots, X_n) \to 0$ in $P_0^n$-probability. Assume that $\mathcal{P}_n$ can be partitioned as $\bigcup_{j=-\infty}^{\infty} \mathcal{P}_{n,j}$ such that, for a sequence $\varepsilon_n \to 0$ with $n\varepsilon_n^2 \to \infty$,*

$$(4.1) \qquad \sum_{j=-\infty}^{\infty} \sqrt{N(\varepsilon_n, \mathcal{P}_{n,j}, d)} \sqrt{\Pi_n(\mathcal{P}_{n,j})} e^{-n\varepsilon_n^2} \to 0,$$

$$(4.2) \qquad \Pi_n(p : P_0 \log(p_0/p) \le \varepsilon_n^2, P_0 \log^2(p_0/p) \le \varepsilon_n^2) \ge e^{-n\varepsilon_n^2}.$$

*Then $\Pi_n(p \in \mathcal{P} : d(p_0, p) > 8\varepsilon_n | X_1, \ldots, X_n) \to 0$ in $P_0^n$-probability.*

Theorem 5 contains a standard posterior convergence theorem (cf. Theorem 2.1 of [8]) as a special case where $\mathcal{P}_n$ is not decomposed (i.e., $\mathcal{P}_{n,0} = \mathcal{P}_n$ and $\mathcal{P}_{n,j} = \varnothing$ for $j \ne 0$), so that $\log N(\varepsilon_n, \mathcal{P}_n, d)$ needs to be bounded by a small multiple of $n\varepsilon_n^2$ in order to satisfy (4.1). At the other extreme, if we decompose $\mathcal{P}_n$ sufficiently finely so that each $\mathcal{P}_{n,j}$ has diameter less than $\varepsilon_n$, then the covering numbers appearing in (4.1) are all 1 and hence (4.1) reduces to

$$(4.3) \qquad \sum_{j=-\infty}^{\infty} \sqrt{\Pi_n(\mathcal{P}_{n,j})} e^{-n\varepsilon_n^2} \to 0.$$

The trade-off between entropy and summability of the square roots of prior probabilities is interesting and requires further investigation; see [15] for a consistency result based on the summability condition.

To prove Theorem 5, we need two auxiliary results. Ghosal, Ghosh and van der Vaart [8] used this result with $\alpha = \beta = 1$.

LEMMA 1. *For any convex set $\mathcal{Q}$ of probability measures with $\inf\{h(P_0, Q): Q \in \mathcal{Q}\} \geq \varepsilon$, any $\alpha, \beta > 0$ and all $n \geq 1$, there exists a test $\phi_n = \phi_n(X_1, \ldots, X_n)$ such that*

$$\sup_{Q \in \mathcal{Q}} (\alpha P_0^n \phi_n + \beta Q^n (1 - \phi_n)) \leq \sqrt{\alpha\beta} e^{-n\varepsilon^2/2}.$$

PROOF. The proof follows by a minor adaptation of a result in [11], pages 475–479, as in [10].  □

COROLLARY 1. *For any set of probability measures $\mathcal{Q}$ with $\inf\{d(P_0, Q): Q \in \mathcal{Q}\} \geq 4\varepsilon$, any $\alpha, \beta > 0$ and all $n \geq 1$, there exists a test $\phi_n$ such that*

$$P_0^n \phi_n \leq \sqrt{\frac{\beta}{\alpha}} N(\varepsilon, \mathcal{Q}, d) \frac{e^{-n\varepsilon^2}}{1 - e^{-n\varepsilon^2}}, \qquad \sup_{Q \in \mathcal{Q}} Q^n (1 - \phi_n) \leq \sqrt{\frac{\alpha}{\beta}} e^{-n\varepsilon^2}.$$

PROOF. For a given $j \in \mathbb{N}$, choose a maximal $j\varepsilon/2$-separated set of points in $S_j = \{Q \in \mathcal{Q}: j\varepsilon < d(Q, P_0) \leq (j+1)\varepsilon\}$. This yields a set $S_j'$ such that the union of the balls of radius $j\varepsilon/2$ centered at these points covers $S_j$. Any such ball $B$ is convex by assumption and satisfies $h(Q, P_0) \geq d(Q, P_0) \geq j\varepsilon/2$ for all $Q \in B$. Because any given ball of radius $j\varepsilon/4$ can contain at most one point of $S_j'$, it follows that $\#S_j' \leq N(\varepsilon j/4, S_j, d) \leq N(\varepsilon, \mathcal{Q}, d)$ for $j \geq 4$. (If $S_j$ is empty, take $S_j'$ empty and adapt the following in the obvious way.)

For every $P_1 \in S_j'$, there exists a test $\omega$ with properties as in Lemma 1, with $\mathcal{Q}$ equal to the ball of radius $j\varepsilon/2$ centered at $P_1$. Let $\phi_n$ be the maximum of all tests attached in this way to some point $P_1 \in S_j'$ for some $j \geq 4$. Then for all $j \geq 4$,

$$P_0^n \phi_n \leq \sum_{j=4}^{\infty} \sum_{P_1 \in S_j'} \sqrt{\frac{\beta}{\alpha}} e^{-nj^2\varepsilon^2/8} \leq \sqrt{\frac{\beta}{\alpha}} N(\varepsilon, \mathcal{Q}, d) \frac{e^{-2n\varepsilon^2}}{1 - e^{-n\varepsilon^2}},$$

$$\sup_{Q \in \bigcup_{i \geq j} S_i} Q^n (1 - \phi_n) \leq \sup_{i \geq j} \sqrt{\frac{\alpha}{\beta}} e^{-ni^2\varepsilon^2/8} \leq \sqrt{\frac{\alpha}{\beta}} e^{-2n\varepsilon^2}.$$  □

PROOF OF THEOREM 5. Clearly, we may assume that the prior charges only $\mathcal{P}_n$. The event $A_n$ that $\int \prod_{i=1}^{n} (p(X_i)/p_0(X_i)) \, d\Pi_n(p) \geq e^{-3n\varepsilon_n^2}$ satisfies $P_0^n(A_n) \to 1$ by Lemma 8.1 of [8] and assumption (4.2). Now, for arbitrary tests $\phi_{n,j}$, we have

$$P_0^n[\Pi_n(P \in \mathcal{P}_{n,j}: d(P, P_0) \geq 8\varepsilon_n | X_1, \ldots, X_n) 1_{A_n}]$$

$$\leq P_0^n \phi_{n,j} + P_0^n \left( (1 - \phi_{n,j}) \int_{\{P \in \mathcal{P}_{n,j} : d(P, P_0) \geq 8\varepsilon_n\}} \prod_{i=1}^n \frac{p(X_i)}{p_0(X_i)} d\Pi_n(p) \right) e^{3n\varepsilon_n^2}$$

$$\leq P_0^n \phi_{n,j} + \sup_{P \in \mathcal{P}_{n,j} : d(P, P_0) \geq 8\varepsilon_n} P^n(1 - \phi_{n,j}) \Pi_n(\mathcal{P}_{n,j}) e^{3n\varepsilon_n^2},$$

which can be bounded by a multiple of

$$\sqrt{\frac{\beta_j}{\alpha_j}} N(2\varepsilon_n, \mathcal{P}_{n,j}, d) e^{-4n\varepsilon_n^2} + \sqrt{\frac{\alpha_j}{\beta_j}} e^{-4n\varepsilon_n^2} \Pi_n(\mathcal{P}_{n,j}) e^{3n\varepsilon_n^2}$$

for the choice of $\phi_{n,j}$ obtained from Corollary 1 with $\varepsilon = 2\varepsilon_n$, $\mathcal{Q} = \{P \in \mathcal{P}_{n,j} : d(P_0, P) \geq 8\varepsilon_n\}$ and any $\alpha_j, \beta_j > 0$. Put $\alpha_j = N(2\varepsilon_n, \mathcal{P}_{n,j}, d)$, $\beta_j = \Pi_n(\mathcal{P}_{n,j})$ and sum over $j$ to obtain the result in view of (4.1). □

## 5. Entropy estimates.

In this section we estimate the entropy of normal mixtures, paying special attention to components with small variances. The main idea is to approximate general normal mixtures by finite mixtures with a small number of components. This same device will also be used in Section 9 to estimate the prior probabilities of Kullback–Leibler type balls and is isolated in the following lemma. The lemma is based on the corresponding one for the super-smooth case (cf. Lemma 3.1 of [9]) and a partitioning argument.

LEMMA 2. *Let $0 < \varepsilon < 1/2$ and $a, \sigma > 0$ be given. For any probability measure $F$ on $[-a, a]$, there exists a discrete probability measure $F'$ on $[-a, a]$ with fewer than $D(a\sigma^{-1} \vee 1) \log \varepsilon^{-1}$ support points, where $D$ is a universal constant, such that*

$$\|p_{F,\sigma} - p_{F',\sigma}\|_\infty \lesssim \frac{\varepsilon}{\sigma}, \qquad \|p_{F,\sigma} - p_{F',\sigma}\|_1 \lesssim \varepsilon (\log \varepsilon^{-1})^{1/2}.$$

PROOF. We can partition the interval $[-a, a]$ into $k = \lfloor 2\alpha/\sigma \rfloor$ disjoint, consecutive subintervals $I_1, \ldots, I_k$ of length $\sigma$ and a final interval $I_{k+1}$ of length $l_{k+1}$ smaller than $\sigma$. We may write $F = \sum_{i=1}^{k+1} F(I_i) F_i$, where each $F_i$ is a probability measure concentrated on $I_i$, and then $p_{F,\sigma} = \sum_{i=1}^{k+1} F(I_i) p_{F_i,\sigma}$. For ease of notation, let $Z_i$ be a random variable distributed according to $F_i$, and for $a_i$ the left endpoint of $I_i$, let $G_i$ be the law of $W_i = (Z_i - a_i)/\sigma$. Thus, $G_i$ is a probability measure on $[0, 1]$ for $i = 1, \ldots, k$ and on $[0, l_{k+1}/\sigma] \subset [0, 1]$ for $i = k + 1$.

By Lemma 3.1 of [9], there exist probability measures $G'_i$ on $[0, 1]$ with $N_i \lesssim \log \varepsilon^{-1}$ support points such that $\|p_{G_i,1} - p_{G'_i,1}\|_\infty \lesssim \varepsilon$. For $i = k + 1$, the measure $G'_i$ can be taken to be supported on $[0, l_{k+1}/\sigma]$. Lemma 3.3 from

the same paper then shows that $\|p_{G_i,1} - p_{G'_i,1}\|_1 \lesssim \varepsilon(\log \varepsilon^{-1})^{1/2}$. Let $F'_i$ be the law of $a_i + \sigma W'_i$ if $W'_i$ has law $G'_i$ and set $F' = \sum_{i=1}^{k+1} F(I_i)F'_i$. Because

$$p_{F_i,\sigma}(x) = \mathrm{E}\phi_\sigma(x - Z_i) = \sigma^{-1}\mathrm{E}\phi((x - a_i)/\sigma - W_i) = \sigma^{-1}p_{G_i,1}((x - a_i)/\sigma),$$

and similarly for $F'_i$ and $G'_i$, we have

$$\|p_{F_i,\sigma} - p_{F'_i,\sigma}\|_\infty = \sigma^{-1}\|p_{G_i,1} - p_{G'_i,1}\|_\infty,$$

$$\|p_{F_i,\sigma} - p_{F'_i,\sigma}\|_1 = \|p_{G_i,1} - p_{G'_i,1}\|_1.$$

Combined with $\|p_{F,\sigma} - p_{F',\sigma}\| \le \sum_{i=1}^{k+1} F(I_i)\|p_{F_i,\sigma} - p_{F'_i,\sigma}\|$, this shows that $p_{F',\sigma}$ has the required distances to $p_{F,\sigma}$. The number of support points of $F'$ is bounded by the number of intervals $k + 1$ times the maximum number of support points of an $F'_i$, and hence is bounded by a multiple of $(a\sigma^{-1} \vee 1)\log \varepsilon^{-1}$.  □

For given numbers $a, b_1, b_2$, let

(5.1)           $$\mathcal{P}_{a,b_1,b_2} = \{p_{F,\sigma} : F[-a,a] = 1, b_1 \le \sigma \le b_2\}.$$

LEMMA 3.   *Let $0 < b_1 < b_2$ and $a > 0$ and define $\mathcal{P}_{a,b_1,b_2}$ by (5.1). Then for $0 < \varepsilon < 1/2$ and $d$ equal to the $L_1$-norm,*

$$\log N(\varepsilon, \mathcal{P}_{a,\tau}, d) \lesssim \log\left(\frac{b_2}{b_1\varepsilon}\right) + \left(\frac{a}{b_1} + 1\right)\left(\log\frac{1}{\varepsilon}\right)\left(\log\frac{1}{\varepsilon} + \log\left(\frac{a}{b_1} + 1\right)\right).$$

*For $d = \|\cdot\|_\infty$, the same bound holds with $\varepsilon$ replaced by $\varepsilon b_1$ if $\varepsilon b_1 < 1/2$. For $d = h$, the same bound holds with $\varepsilon$ replaced by $\varepsilon^2$.*

PROOF.   Let $\alpha, \beta < 1/2, \gamma$ and $\delta$ be given positive numbers. Fix a minimal $\alpha$-net $\Sigma$ over the interval $[b_1, b_2]$ and let $\mathcal{F}$ be the set of discrete probability distributions on $[-a, a]$ with at most $N \le D(ab_1^{-1} \vee 1)\log \beta^{-1}$ support points, for the constant $D$ of Lemma 2. For every $\sigma \in [b_1, b_2]$, there exists $\sigma' \in \Sigma$ with $|\sigma - \sigma'| \le \alpha$, whence

$$\|p_{F,\sigma} - p_{F,\sigma'}\|_\infty \le \|\phi_\sigma - \phi_{\sigma'}\|_\infty \lesssim |\sigma - \sigma'|/(\sigma \wedge \sigma')^2 \lesssim \alpha/b_1^2,$$

$$\|p_{F,\sigma} - p_{F,\sigma'}\|_1 \le \|\phi_\sigma - \phi_{\sigma'}\|_1 \lesssim |\sigma - \sigma'|/(\sigma \wedge \sigma') \lesssim \alpha/b_1.$$

By Lemma 2, for a sufficiently large $D$, there exists, for every given probability measure $F$ on $[-a, a]$, an element $F' \in \mathcal{F}$ (possibly depending on $\sigma'$) such that $\|p_{F,\sigma'} - p_{F',\sigma'}\|_\infty \lesssim \beta/b_1$ and $\|p_{F,\sigma'} - p_{F',\sigma'}\|_1 \lesssim \beta(\log \beta^{-1})^{1/2}$. Hence,

$$\|p_{F,\sigma} - p_{F',\sigma'}\|_\infty \lesssim \alpha/b_1^2 + \beta/b_1,$$

$$\|p_{F,\sigma} - p_{F',\sigma'}\|_1 \lesssim \alpha/b_1 + \beta(\log \beta^{-1})^{1/2}.$$

Thus, $\mathcal{P} = \{p_{F,\sigma} : (F,\sigma) \in \mathcal{F} \times \Sigma\}$ is an $\varepsilon$-net over $\mathcal{P}_{a,b_1,b_2}$ for $\|\cdot\|_\infty$ and $\|\cdot\|_1$, respectively, if the expressions above are made to be less than $\varepsilon$.

We next construct a finite net over $\mathcal{P}$ by restricting the support points and weights of $F$ to suitable grids. For a fixed $\gamma$-net $\mathcal{S}$ over the $N$-dimensional simplex for the $\ell_1$-norm, let $\mathcal{P}'$ be the set of all $p_{F,\sigma} \in \mathcal{P}$ such that the $N$ support points of $F$ are among the points $\{0, \pm\delta, \pm2\delta, \dots\} \cap [-a,a]$ with weights belonging to $\mathcal{S}$. We may project $p_{F,\sigma} \in \mathcal{P}$ into $p_{F',\sigma} \in \mathcal{P}'$ by first moving the point masses of $F$ to a closest point in the grid $\{0, \pm\delta, \pm2\delta, \dots\}$ and then changing the vector of sizes of the point masses to a closest vector in $\mathcal{S}$. Let $z_1, z_2, \dots, z_1', z_2', \dots, p_1, p_2, \dots, p_1', p_2', \dots$ be such that $F = \sum p_j \delta_{z_j}$ and $F' = \sum p_j \delta_{z_j'}$. Then

$$|p_{F,\sigma}(x) - p_{F',\sigma}(x)|$$
$$\leq \sum_j \{p_j |\phi_\sigma(x-z_j) - \phi_\sigma(x-z_j')| + |p_j - p_j'| \phi_\sigma(x-z_j')\}.$$

Thus, $\|p_{F,\sigma} - p_{F',\sigma}\|_\infty \lesssim \delta \|\phi_\sigma'\|_\infty + \gamma \|\phi_\sigma\|_\infty \lesssim \delta/b_1^2 + \gamma/b_1$ and $\|p_{F,\sigma} - p_{F',\sigma}\|_1 \lesssim \delta/b_1 + \gamma$. Hence, $\mathcal{P}'$ is a $c\eta$-net over $\mathcal{P}_{a,b_1,b_2}$ for $c$ a universal constant and for $\eta = \eta_\infty = \alpha/b_1^2 + \beta/b_1 + \delta/b_1^2 + \gamma/b_1$ for $\|\cdot\|_\infty$ and $\eta = \eta_1 = \alpha/b_1 + \beta(\log \beta^{-1})^{1/2} + \delta/b_1 + \gamma$ for $\|\cdot\|_1$. Because $h^2(f,g) \leq \|f-g\|_1$ for any two densities, $\mathcal{P}'$ is a Hellinger $c^2\eta^2$-net over $\mathcal{P}_{a,b_1,b_2}$ for $\eta = \eta_1$.

There are at most $((b_2 - b_1)/\alpha) \vee 1$ possible choices of $\sigma \in \Sigma$. Each $z_j'$ can assume at most $2a/\delta + 1$ different values, $j = 1, \dots, N$. The cardinality of a minimal $\gamma$-net $\mathcal{S}$ over the $N$-dimensional unit simplex is bounded by $(5/\gamma)^N$ for $\gamma \leq 1$ (cf. Lemma A.4 of [9]). Therefore,

$$\#\mathcal{P}' \lesssim \left(\frac{b_2 - b_1}{\alpha} \vee 1\right) \times \left(\frac{2a}{\delta} + 1\right)^N \times \left(\frac{5}{\gamma \wedge 1}\right)^N.$$

Because $N \leq Dab_1^{-1} \log \beta^{-1}$, by construction, it follows that

$$\log(\#\mathcal{P}') \lesssim \log\left(\frac{b_2 - b_1}{\alpha} \vee 1\right) + \left(\frac{a}{b_1} \vee 1\right) \log \frac{1}{\beta} \left[\log\left(\frac{2a}{\delta} + 1\right) + \log\left(\frac{5}{\gamma \wedge 1}\right)\right].$$

This number, with $\alpha = \delta = b_1 \varepsilon$, $\beta = \varepsilon (\log \varepsilon^{-1})^{-1/2}$ and $\gamma = \varepsilon$ for given $\varepsilon < 1/2$, is a bound on the $D'\varepsilon$-entropy for $\|\cdot\|_1$ for a universal constant $D'$, and with $\alpha = \delta = b_1^2 \varepsilon$, $\beta = b_1 \varepsilon$ and $\gamma = b_1 \varepsilon$, a bound on the $D'\varepsilon$-entropy for $\|\cdot\|_\infty$ is obtained upon simplification. (If $D' > 1$, we replace $\varepsilon$ by $\varepsilon/D'$.)  □

For positive numbers $a, \tau, b_1 < b_2$, let

(5.2)        $$\mathcal{P}_{a,\tau} = \{p_{F,\sigma} : F[-a,a] = 1, b_1 \tau \leq \sigma \leq b_2 \tau\},$$

where $\tau$ is small.

THEOREM 6. *Let $b_1 < b_2$, $\tau < 1/4$ and $a \geq e$ be given positive numbers and define $\mathcal{P}_{a,\tau}$ by (5.2). Then, for $0 < \varepsilon < 1/2$ and $d$ the $L_1$-norm or Hellinger distance,*

$$\log N(\varepsilon, \mathcal{P}_{a,\tau}, d) \lesssim \frac{a}{\tau}\left(\log\frac{1}{\varepsilon}\right)\left(\log\frac{a}{\varepsilon\tau}\right), \tag{5.3}$$

*where the constant in "$\lesssim$" depends on $b_1, b_2$ only. For $d = \|\cdot\|_\infty$, (5.3) holds with $\log\varepsilon^{-1}$ replaced by $\log(\varepsilon\tau)^{-1}$. Further, for any of the three metrics,*

$$\log N_{[\cdot]}(\varepsilon, \mathcal{P}_{a,\tau}, d) \lesssim \frac{a}{\tau}\left(\log\frac{a}{\varepsilon\tau}\right)^2. \tag{5.4}$$

PROOF. Inequality (5.3) can be deduced from Lemma 3 by replacing $b_1$ and $b_2$ by $b_1\tau$ and $b_2\tau$ and then simplifying the resulting entropy bounds.

To obtain the bound on the bracketing numbers of $\mathcal{P}_{a,\tau}$, we first note that $H(x) = (b_1\tau)^{-1}\phi(x/(2b_2\tau))\mathbf{1}\{|x| > 2a\} + (b_1\tau)^{-1}\phi(0)\mathbf{1}\{|x| \leq 2a\}$ is an envelope for $\mathcal{P}_{a,\tau}$. Given an $\eta$-net $f_1, \ldots, f_M$ for $\|\cdot\|_\infty$, the brackets $[l_i, u_i]$, where $l_i = (f_i - \eta) \vee 0$ and $u_i = (f_i + \eta) \wedge H$, cover $\mathcal{P}_{a,\tau}$. Thus, $u_i - l_i \leq (2\eta) \wedge H$ and the size of these brackets in $L_1$ can be bounded by

$$\int (u_i - l_i)\, d\lambda \lesssim \|u_i - l_i\|_\infty B + \int_{|x| > B} H(x)\, dx \lesssim \eta B + \phi(B/2b_2\tau),$$

for any $B > 2a$, by the tail bound for the normal distribution. For $B = (b_2 \vee 1)2a(\log\eta^{-1})^{1/2}$, we obtain the upper bound equal to a multiple of $\eta(\log\eta^{-1})^{1/2}a + \eta^{8a^2} \lesssim \eta(\log\eta^{-1})^{1/2}a$. Thus, there exists a constant $D$ (possibly depending on $b_1, b_2$) such that the $D\eta(\log\eta^{-1})^{1/2}a$-bracketing number for the $L_1$-norm is bounded by the uniform $\eta$-covering number obtained previously. Choose $\eta = D\varepsilon a^{-2}(\log\varepsilon^{-1})^{-1/2}$ for an appropriate constant $D$ and simplify to obtain (5.4)   □

The main difference between the bound given by Theorem 6 and the bound when the scale is bounded away from zero (cf. [9]) is the presence of the factor $a\tau^{-1}$. This factor is the main driving force for the slower rate of convergence of the posterior in the present situation compared to the super-smooth case.

The set of mixtures with an arbitrary mixing distribution on $\mathbb{R}$ is not totally bounded and hence Theorem 6 can only be extended to mixing distributions with possibly noncompact support if the mixing distribution is restricted in some other way. We shall extend the theorem to mixtures with mixing distributions whose tails are bounded by a given function (such as the normal density).

For a given decreasing function $A:(0,\infty) \to [0,1]$ with inverse $A^{-1}$ and positive numbers $\tau$ and $b_1 < b_2$, we consider the class of densities

$$\mathcal{P}_{A,\tau} = \{p_{F,\sigma} : F[-a,a]^c \leq A(a) \text{ for all } a, b_1\tau \leq \sigma \leq b_2\tau\}. \tag{5.5}$$

For $\Phi$ the standard normal distribution function, let $\bar{A}$ be the function defined by $\bar{A}(a) = A(a) + \int_{a/2}^{\infty} A\, d\lambda + 1 - \Phi(a)$.

**Theorem 7.** *Let $b_1 < b_2$ and $\tau < 1/4$ be given positive numbers, let $A : (0, \infty) \to [0, 1]$ be a decreasing function and define $\mathcal{P}_{A,\tau}$ as in (5.5). Then, for $0 < \varepsilon < \min(1/4, A(e))$, we have that*

$$\log N(3\varepsilon, \mathcal{P}_{A,\tau}, \|\cdot\|_1) \lesssim \frac{A^{-1}(\varepsilon)}{\tau} \left( \log \frac{1}{\varepsilon} \right) \left( \log \frac{A^{-1}(\varepsilon)}{\varepsilon\tau} \right),$$

*where the constant in "$\lesssim$" depends on $b_1, b_2$ only. Furthermore, for a constant $c$ depending on $b_1, b_2$ only,*

$$\log N_{[\cdot]}(c\varepsilon, \mathcal{P}_{A,\tau}, \|\cdot\|_1) \lesssim \frac{\bar{A}^{-1}(\varepsilon\tau)}{\tau} \left( \log \frac{\bar{A}^{-1}(\varepsilon\tau)}{\varepsilon\tau} \right)^2.$$

*For the entropy relative to $h$, the same bounds hold with $\varepsilon$ replaced by $\varepsilon^2$.*

PROOF. Because $a_\varepsilon = A^{-1}(\varepsilon)$ satisfies $F[-a_\varepsilon, a_\varepsilon]^c \le \varepsilon$ for every $F$ as in the definition of $\mathcal{P}_{A,\tau}$, Lemma A.3 of [9] shows that the $L_1$-distance between $\mathcal{P}_{A,\tau}$ and $\mathcal{P}_{a_\varepsilon,\tau}$ is bounded above by $2\varepsilon$. It follows that an $\varepsilon$-net over $\mathcal{P}_{a_\varepsilon,\tau}$ is a $3\varepsilon$-net over $\mathcal{P}_{A,\tau}$. In view of Theorem 6, this implies the bound on the entropy without bracketing given in the first inequality of the theorem.

To bound the bracketing numbers, we obtain by partial integration, for any $x > a > 0$,

$$\int_a^{\infty} \phi_\sigma(x - z)\, dF(z) = (1 - F)(a)\phi_\sigma(x - a) + \int_a^{\infty} \phi_\sigma'(z - x)(1 - F)(z)\, dz$$

$$\lesssim A(a)\phi_\sigma(x - a) + \frac{\phi(-x)}{\sigma} + \frac{A(x/2)}{\sigma}.$$

For $x < -a < 0$, the same bound is valid for $\int_{-\infty}^{-a} \phi_\sigma(x - z)\, dF(z)$, but with $\phi_\sigma(x - a)$ replaced by $\phi_\sigma(x + a)$ and with $-x$ replaced by $x$. Also, for $a > 0$ and $x < a$,

$$\int_a^{\infty} \phi_\sigma(x - z)\, dF(z) \le \phi_\sigma(x - a)F[a, \infty) \lesssim \phi_\sigma(x - a)A(a).$$

For $x > -a$, the same bound is valid for $\int_{-\infty}^{-a} \phi_\sigma(x - z)\, dF(z)$, but with $\phi_\sigma(x - a)$ replaced by $\phi_\sigma(x + a)$. Hence, for any $x$ and $a > 0$,

$$\int_{|z| > a} \phi_\sigma(x - z)\, dF(z) \lesssim (\phi_{b_2\tau}(x - a) + \phi_{b_2\tau}(x + a))A(a)$$

$$+ \frac{1}{\tau}(\phi(x) + A(x/2))\mathbf{1}\{|x| \ge a\}.$$

Let $H$ denote the appropriate multiple of the right-hand side of this inequality. If $F_a$ is the renormalized restriction to $[-a, a]$ of a probability measure $F$, then $F[-a, a]p_{F_a,\sigma} \le p_{F,\sigma} \le F[-a, a]p_{F_a,\sigma} + H$. Consequently, given $\varepsilon$-brackets $[l_i, u_i]$ that cover $\mathcal{P}_{a,\tau}$, there exists for every $(F, \sigma)$, as in the definition of $\mathcal{P}_{A,\tau}$, a bracket $[l_i, u_i]$ with $l_i(1 - A(a)) \le p_{F,\sigma} \le F[-a, a]u_i + H \le u_i + H$. Thus, the brackets $[l_i(1 - A(a)), u_i + H]$ cover $\mathcal{P}_{A,\tau}$. The size in $L_1$ of a bracket $[l_i, u_i]$ is bounded by $\|u_i - l_i\|_1 + \|l_i\|_1 A(a) + \|H\|_1 \lesssim \varepsilon + \bar{A}(a)/\tau$. We now choose $a$ such that $\bar{A}(a) \le \tau\varepsilon$ and apply Theorem 6 to bound the number of brackets $[l_i, u_i]$. □

As an example, if the mixing distributions have sub-Gaussian tails, then we can apply the preceding theorem with $A$ equal to $1 - \Phi(a) \le \phi(a)$, whence $\bar{A}$ is bounded by a multiple of the same function. Then, both $A^{-1}(\varepsilon)$ and $\bar{A}^{-1}(\varepsilon)$ are bounded by a multiple of $(\log \varepsilon^{-1})^{1/2}$ and the (bracketing) entropy is bounded by a multiple of $\tau^{-1}(\log(\varepsilon\tau)^{-1})^{5/2}$. Provided the tails of the mixing distributions are bounded by a function of the form $A(a) = e^{-da^\delta}$, the entropy of the set of mixtures increases at most through a power of $\log(\varepsilon\tau)^{-1}$. On the other hand, polynomially decreasing tails incur an additional factor of $\varepsilon^{-k}$ in the entropy bounds of Theorem 7.

**6. Approximation results.** If $p_0$ is a twice differentiable density, then for $d$ equal to the $L_p$-distance, it is well known that $d(p_0, p_0 * \phi_\sigma) = O(\sigma^2)$. In the following lemma, we establish this for $d = h$.

LEMMA 4.  *Let $p_0$ be a twice continuously differentiable probability density.*

- *If $\int (p_0''/p_0)^2 p_0 \, d\lambda < \infty$ and $\int (p_0'/p_0)^4 p_0 \, d\lambda < \infty$, then $h(p_0, p_0 * \phi_\sigma) \lesssim \sigma^2$.*
- *If $p_0$ is bounded with $\int |p_0''| \, d\lambda < \infty$, then $\|p_0 - p_0 * \phi_\sigma\|_1 \lesssim \sigma^2$.*

*In both cases, the constants in "$\lesssim$" depend on the given integrals only.*

PROOF.   By the assumption of $P_0$-integrability of the functions $p_0'/p_0$ and $p_0''/p_0$, we have $\int |p_0^{(i)}| \, d\lambda < \infty$ for $i = 1, 2$. Therefore, $p_0$ and $p_0'$ are uniformly bounded, from which it can be seen that $p_\sigma(x) = \int p_0(x - \sigma y)\phi(y) \, dy$ is twice partially differentiable relative to $\sigma$, with derivatives $\dot{p}_\sigma(x)$ and $\ddot{p}_\sigma(x)$ given by

$$\dot{p}_\sigma(x) = -\int p_0'(x - \sigma y)y\phi(y) \, dy,$$

$$\ddot{p}_\sigma(x) = -\int p_0''(x - \sigma y)y^2\phi(y) \, dy.$$

Using Taylor's theorem with the integral form of the remainder (cf. [2], page 120), we have

$$p_\sigma^{1/2}(x) - p_0^{1/2}(x) = \sigma \frac{\dot{p}_0(x)}{2p_0^{1/2}(x)} + \frac{1}{2}\sigma^2 \int_0^1 \left( \frac{\ddot{p}_{s\sigma}(x)}{p_{s\sigma}^{1/2}(x)} - \frac{1}{2} \frac{\dot{p}_{s\sigma}^2(x)}{p_{s\sigma}^{3/2}(x)} \right)(1-s)\,ds.$$

Because $\dot{p}_0(x) = -\int p_0'(x)y\phi(y)\,dy = 0$ for every $x$, we obtain

$$h^2(p_\sigma, p_0) = \frac{1}{4}\sigma^4 \int \left( \int_0^1 \left( \frac{\ddot{p}_{s\sigma}(x)}{p_{s\sigma}^{1/2}(x)} - \frac{1}{2} \frac{\dot{p}_{s\sigma}^2(x)}{p_{s\sigma}^{3/2}(x)} \right)(1-s)\,ds \right)^2 dx$$

$$\leq \frac{1}{2}\sigma^4 \int_0^1 \int \left[ \left( \frac{\ddot{p}_{s\sigma}(x)}{p_{s\sigma}^{1/2}(x)} \right)^2 + \frac{1}{4}\left( \frac{\dot{p}_{s\sigma}^2(x)}{p_{s\sigma}^{3/2}(x)} \right)^2 \right] dx \times (1-s)^2\,ds.$$

Now, for any $\sigma$, by the Cauchy–Schwarz inequality,

$$\ddot{p}_\sigma^2(x) = \left( \int \frac{p_0''(x-\sigma y)}{p_0^{1/2}(x-\sigma y)} y^2 p_0^{1/2}(x-\sigma y)\phi(y)\,dy \right)^2$$

$$\leq \int \frac{(p_0''(x-\sigma y))^2}{p_0(x-\sigma y)} y^4\,dy \times p_\sigma(x).$$

Furthermore, by Hölder's inequality with $p=4$ and $q=4/3$, we have

$$\dot{p}_\sigma^4(x) = \left( \int \frac{p_0'(x-\sigma y)}{p_0^{3/4}(x-\sigma y)} y p^{3/4}(x-\sigma y)\phi(y)\,dy \right)^2$$

$$\leq \int \left( \frac{p_0'(x-\sigma y)}{p_0^{3/4}(x-\sigma y)} \right)^4 y^4 \phi(y)\,dy \times \left( \int p(x-\sigma y)\phi(y)\,dy \right)^3$$

$$= \int \left( \frac{p_0'(x-\sigma y)}{p_0^{3/4}(x-\sigma y)} \right)^4 y^4\,dy \times p_\sigma^3(x).$$

The required bound for the proof of the first assertion now follows by substituting these inequalities into the expression for the Hellinger distance and interchanging the order of integration.

The proof of the second assertion is similar, but easier. $\quad\square$

The following lemma bounds the distance between a normal mixture with a mixing distribution $F$ and that with a discrete approximation to $F$. This result, which extends Lemma 5.1 of [9], will be instrumental in lower-bounding the prior probability of Kullback–Leibler-type neighborhoods.

LEMMA 5. *Let* $\mathbb{R} = \bigcup_{j=0}^N U_j$ *be a partition of* $\mathbb{R}$ *and* $F' = \sum_{j=1}^N p_j \delta_{z_j}$ *be a probability measure with* $z_j \in U_j$ *for* $j = 1, \ldots, N$. *Then, for any probability* *measure $F$ on $\mathbb{R}$,*

$$\|p_{F,\sigma} - p_{F',\sigma}\|_\infty \lesssim \frac{1}{\sigma^2} \max_{1 \le j \le N} \lambda(U_j) + \frac{1}{\sigma} \sum_{j=1}^N |F(U_j) - p_j|,$$

$$\|p_{F,\sigma} - p_{F',\sigma}\|_1 \lesssim \frac{1}{\sigma} \max_{1 \le j \le N} \lambda(U_j) + \sum_{j=1}^N |F(U_j) - p_j|.$$

PROOF. Bound $p_{F,\sigma}(x) - p_{F',\sigma}(x)$ by

$$\int_{U_0} \phi_\sigma(x - z) \, dF(z) + \sum_{j=1}^N \int_{U_j} (\phi_\sigma(x - z) - \phi_\sigma(x - z_j)) \, dF(z)$$

$$+ \sum_{j=1}^N \phi_\sigma(x - z_j)(F(U_j) - p_j).$$

The result now follows because $F(U_0) = 1 - \sum_{j=1}^N F(U_j) \le \sum_{j=1}^N |F(U_j) - p_j|$, $\|\phi_\sigma\|_\infty \lesssim \sigma^{-1}$, $\|\phi'_\sigma\|_\infty \lesssim \sigma^{-2}$ and $\|\phi_\sigma(\cdot - z) - \phi_\sigma(\cdot - z')\|_1 \lesssim \sigma^{-1}|z - z'|$. □

Bounds on the Kullback–Leibler divergence require some control on quotients of the type $p_0/p_{F,\sigma}$. The following lemma, which is implicit in Remark 3 of [7], is useful for this purpose.

LEMMA 6. *Let $p$ be a bounded probability density such that $p$ is nondecreasing on $(-\infty, a]$, bounded away from 0 on $[a, b]$ and nonincreasing on $[b, \infty)$ for some $a \le b$. Then, for every $\tau > 0$, there exists a constant $C > 0$ such that $p * \phi_\sigma \ge Cp$ for every $\sigma < \tau$.*

The next three lemmas are useful to control the Kullback–Leibler divergence and similar quantities in terms of the Hellinger distance. The first lemma is a simplification of Theorem 5 of [16].

LEMMA 7. *For every $b > 0$, there exists a constant $\varepsilon_b > 0$ such that for all probability measures $P$ and $Q$ with $0 < h^2(p,q) < \varepsilon_b P(p/q)^b$, with $\log_+ x = \log x \vee 0$,*

$$P \log \frac{p}{q} \lesssim h^2(p,q) \left\{ 1 + \frac{1}{b} \log_+ \frac{1}{h(p,q)} + \frac{1}{b} \log_+ P\left(\frac{p}{q}\right)^b \right\},$$

$$P \left( \log \frac{p}{q} \right)^2 \lesssim h^2(p,q) \left\{ 1 + \frac{1}{b} \log_+ \frac{1}{h(p,q)} + \frac{1}{b} \log_+ P\left(\frac{p}{q}\right)^b \right\}^2.$$

PROOF. The function $r:(0, \infty) \to \mathbb{R}$ defined implicitly by $\log x = 2(\sqrt{x} - 1) - r(x)(\sqrt{x} - 1)^2$ possesses the following properties:

(i) $r$ is nonnegative and decreasing.

(ii) $r(x) \sim \log x^{-1}$ as $x \downarrow 0$, whence there exists $\varepsilon' > 0$ such that $r(x) \leq 2\log x^{-1}$ on $[0, \varepsilon']$ (a computer graph indicates that $\varepsilon' = 0.4$ will suffice).

(iii) For every $b > 0$, there exists $\varepsilon_b'' > 0$ such that $x^b r(x)$ is increasing on $[0, \varepsilon_b'']$. (For $b \geq 1$, we may take $\varepsilon_b'' = 1$, but for $b$ close to zero, $\varepsilon_b''$ must be very small.)

Using these properties and $h^2(p, q) = -2P(\sqrt{q/p} - 1)$, we obtain

$$P \log \frac{p}{q} = h^2(p, q) + P\left[ r\left(\frac{q}{p}\right)\left(\sqrt{\frac{q}{p}} - 1\right)^2 \right]$$

$$\leq h^2(p, q) + r(\varepsilon)h^2(p, q) + P\left[ r\left(\frac{q}{p}\right)1\left\{ \frac{q}{p} \leq \varepsilon \right\} \right]$$

$$\leq h^2(p, q) + 2\log\frac{1}{\varepsilon}h^2(p, q) + 2\varepsilon^b \log\frac{1}{\varepsilon}P\left(\frac{p}{q}\right)^b$$

for $\varepsilon \leq \varepsilon' \wedge \varepsilon_b'' \wedge 4$. The proof of the first inequality now follows by choosing $\varepsilon^b = h^2(p, q)/P(p/q)^b$ and $\varepsilon_b = (\varepsilon' \wedge \varepsilon_b'' \wedge 4)^b$.

To prove the second inequality, note that $|\log x| \leq 2|\sqrt{x} - 1|$, $x \geq 1$, and so

$$P\left[ \left( \log\frac{p}{q} \right)^2 1\left\{ \frac{q}{p} \geq 1 \right\} \right] \leq 4P\left( \sqrt{\frac{q}{p}} - 1 \right)^2 = 4h^2(p, q).$$

Next, for $\varepsilon \leq \varepsilon_{b/2}''$, in view of the third property of $r$ we have

$$P\left[ \left( \log\frac{p}{q} \right)^2 1\left\{ \frac{q}{p} \leq 1 \right\} \right] \leq 8P\left( \sqrt{\frac{q}{p}} - 1 \right)^2 + 2P\left[ r^2\left(\frac{q}{p}\right)\left( \sqrt{\frac{q}{p}} - 1 \right)^4 1\left\{ \frac{q}{p} \leq 1 \right\} \right]$$

$$\leq 8h^2(p, q) + 2r^2(\varepsilon)h^2(p, q) + 2\varepsilon^b r^2(\varepsilon)P\left(\frac{p}{q}\right)^b.$$

With $\varepsilon^b = h^2(p, q)/P(p/q)^b$ and $\varepsilon_b \leq (\varepsilon' \wedge \varepsilon_{b/2}'')^b$, the proof follows from (ii). $\square$

The next lemma is the limiting case of the preceding lemma as $b \uparrow \infty$. The first assertion was proved by Birgé and Massart ([1], equation (7.6)). The second assertion improves on Lemma 8.3 of [8].

LEMMA 8. *For every pair of probability densities $p$ and $q$,*

$$P \log \frac{p}{q} \lesssim h^2(p, q)\left( 1 + \log\left\| \frac{p}{q} \right\|_\infty \right),$$

$$P\left( \log\frac{p}{q} \right)^2 \lesssim h^2(p, q)\left( 1 + \log\left\| \frac{p}{q} \right\|_\infty \right)^2.$$

PROOF. It can be checked that, for $b \geq 1$, we can choose $\varepsilon_b'' = 1$ in the preceding proof. Furthermore, we can choose $\varepsilon' = 1$ if we use the bound $r(x) \leq 2 + 2\log x$ rather than the bound $r(x) \leq 2\log x$. This leads to the same types of bound as in Lemma 7, which are then seen to be valid for every $b \geq 2$ and any probability densities $p$ and $q$ with $h^2(p,q) \leq P(p/q)^b$, since $\varepsilon_b \geq 1$. Here, $P(p/q)^b = Q(p/q)^{b+1} \geq (Q(p/q))^{b+1} \geq 1$ for $b > 1$ by Jensen's inequality. Thus, the bounds of Lemma 7 hold for every sufficiently large $b$ and every $p$ and $q$ with $h^2(p,q) \leq 1$. For $b \uparrow \infty$, we have that $(P(p/q)^b)^{1/b}$ converges to the $L_\infty(P)$-norm of $p/q$, and the bounds tend to the bounds given by the present lemma. $\square$

Given the control of the supremum of likelihood ratios, we can also compare the Kullback–Leibler divergences of two densities relative to a third density.

LEMMA 9. *For any probability densities $p$, $q$ and $r$,*

$$P \log \frac{p}{r} \leq P \log \frac{p}{q} + 2h(q,r) \left\| \frac{p}{r} \right\|_\infty^{1/2},$$

$$P \left( \log \frac{p}{r} \right)^2 \leq 4P \left( \log \frac{p}{q} \right)^2 + 16h^2(p,q) + 16h^2(q,r) \left\| \frac{p}{r} \right\|_\infty + 16h^2(p,r).$$

*Here, $p/r$ is read as $0$ if $p = 0$ and as $\infty$ if $r = 0 < p$.*

PROOF. To prove the first relation, write $P \log(p/r)$ as the sum of $P \log(p/q)$ and $P \log(q/r)$. Using $\log x \leq 2(\sqrt{x} - 1)$ and the Cauchy–Schwarz inequality, we obtain

$$P \log \frac{q}{r} \leq 2P \frac{1}{\sqrt{r}} (\sqrt{q} - \sqrt{r})$$

$$\leq 2 \left\| \frac{p}{r} \right\|_\infty^{1/2} \int \sqrt{p} (\sqrt{q} - \sqrt{r}) \, d\lambda$$

$$\leq 2 \left\| \frac{p}{r} \right\|_\infty^{1/2} h(q,r).$$

By the relations $\log_+ x \leq 2|\sqrt{x} - 1|$ and $\log_- x = \log_+(1/x) \leq 2|\sqrt{1/x} - 1|$, we have, for any probability densities $p$, $q$, $r$,

$$P \left( \log_+ \frac{q}{r} \right)^2 \leq 4P \left( \sqrt{\frac{q}{r}} - 1 \right)^2 \leq 4 \left\| \frac{p}{r} \right\|_\infty h^2(q,r),$$

$$P \left( \log_- \frac{p}{r} \right)^2 \leq 4P \left( \sqrt{\frac{r}{p}} - 1 \right)^2 = 4h^2(p,r).$$

Since $|\log \frac{p}{r}| \leq \log_+ \frac{p}{q} + \log_- \frac{p}{q} + \log_+ \frac{q}{r} + \log_- \frac{p}{r}$, the second relation follows from the triangle inequality for the $L_2(P)$-norm. $\quad\square$

**7. Prior estimates.** The following extension of Lemma 6.1 of [8] gives a useful probability bound.

LEMMA 10. *For given $N \in \mathbb{N}$, let $(p_1, \ldots, p_N)$ be an arbitrary point in the $N$-dimensional unit simplex and let $(X_1, \ldots, X_N)$ be Dirichlet distributed with parameters $(\alpha_1, \ldots, \alpha_N)$, with $\alpha_j \leq 1$ for every $j$ and $\sum_{j=1}^N \alpha_j = m$. Let $a$ and $b$ be positive numbers. Then, for every $0 < \varepsilon < 1/4$ with $\varepsilon^b \leq a\alpha_j$ and $\varepsilon N \leq 1$, and constants $c$ and $C$ that depend only on $a, b, m$,*

$$(7.1) \qquad \Pr\left(\sum_{j=1}^N |X_j - p_j| \leq 2\varepsilon, \min_{1 \leq j \leq N} X_j \geq \frac{\varepsilon^2}{2}\right) \geq C e^{-cN \log \varepsilon^{-1}}.$$

PROOF. As in the proof of Lemma 6.1 of [8], we can assume without loss of generality that $p_N \geq N^{-1}$, and if $|x_j - p_j| \leq \varepsilon^2$ for $j = 1, \ldots, N-1$, then $x_N \geq \varepsilon^2$ and $\sum_{j=1}^N |x_j - p_j| \leq 2\varepsilon$. Using $\Gamma(\alpha) = \Gamma(1 + \alpha)/\alpha \leq 1/\alpha$ for $0 < \alpha \leq 1$ and the fact that $\alpha_j \geq \varepsilon^b/a$, we obtain

$$\Pr\left(|X_j - p_j| \leq \varepsilon^2, X_j \geq \frac{\varepsilon^2}{2}, j = 1, \ldots, N\right)$$
$$\geq \frac{\Gamma(m)}{\prod_{j=1}^N \Gamma(\alpha_j)} \prod_{j=1}^{N-1} \int_{\max((p_j - \varepsilon^2), \varepsilon^2/2)}^{\min((p_j + \varepsilon^2), 1)} x_j^{\alpha_j - 1} \, dx_j$$
$$\geq \Gamma(m)(\varepsilon^2/2)^{(N-1)} (\varepsilon^b/\alpha)^N$$
$$\geq C \exp(-cN \log \varepsilon^{-1}). \qquad\qquad\qquad\square$$

**8. Tail mass of Dirichlet posterior.** To obtain a posterior convergence rate $\varepsilon_n$ for Dirichlet mixtures, we need to show for some sufficiently large $a$ that

$$(8.1) \qquad \mathrm{E}\Pi_n(F : F[-2a, 2a]^c > \varepsilon_n^2 | X_1, \ldots, X_n) \to 0.$$

In [9], (8.1) was derived by showing that the prior mass of the set in the display is exponentially small. This forces us to increase $a$ with $n$ sufficiently fast, but in the present situation, this method would lead to very restrictive tail conditions on the Dirichlet base measure. Instead, when the true distribution is compactly supported, we shall verify (8.1) for a fixed large $a$ by calculations using explicit properties of the Dirichlet prior and posterior.

LEMMA 11. *Let the true distribution of $X_1, \ldots, X_n$ be i.i.d. $P_0$. If the model is as described in Section 2 and $\alpha$ has a positive and continuous density on $[-a, a]$, then for any $\varepsilon > 0$ and $0 < b < a\sigma_n^{-1}$, there exists $K$ not depending on $n$ such that*

$$\mathrm{E}\left[\Pr(F[-2a, 2a]^c > \varepsilon | X_1, \ldots, X_n) 1\left\{\max_{1 \leq i \leq n} |X_i| \leq a\right\}\right]$$

$$\lesssim \mathrm{E}\Pr(\sigma > b\sigma_n | X_1, \ldots, X_n) + \frac{\alpha[-2a, 2a]^c}{\varepsilon(\alpha(\mathbb{R}) + n)} + Kn\varepsilon^{-1}e^{-a^2/4b^2\sigma_n^2}.$$

*Moreover, if $P_0$ is compactly supported and satisfies the assumptions in Section 1.2, $\alpha$ has positive and continuous density on an interval containing the support of $P_0$, $b_n \to \infty$ is a sequence with $b_n\sigma_n \to 0$, $n\varepsilon_n^{-2}e^{-a^2/4b_n^2\sigma_n^2} \to 0$ and $\Pr(\sigma > b_n\sigma_n) = o(e^{-n\varepsilon_n^2})$ for a sequence $\varepsilon_n$ such that (4.2) holds, then (8.1) holds.*

PROOF. To prove the lemma, it is useful to describe the Dirichlet prior and the observations from the Dirichlet mixtures structurally as follows:

- $F \sim D_\alpha$ and $\sigma/\sigma_n \sim G$, independently;
- given $(F, \sigma)$, the variables $\theta_1, \ldots, \theta_n$ are an i.i.d. sample from $F$;
- given $(F, \sigma, \theta_1, \ldots, \theta_n)$, the variables $e_1, \ldots, e_n$ are i.i.d. $N(0, \sigma^2)$;
- the variables $X_1, \ldots, X_n$ are defined as $X_i = \theta_i + e_i$.

Let $G_n(s) = G(s/\sigma_n)$. Given $(\theta_1, \ldots, \theta_n)$, the observations $X_1, \ldots, X_n$ are independent of $F$ and hence the conditional distribution of $F$ given $(X_1, \ldots, X_n, \theta_1, \ldots, \theta_n)$ is independent of $X_1, \ldots, X_n$. This allows us to write

$$\Pr(F[-2a, 2a]^c > \varepsilon | X_1, \ldots, X_n)$$

$$= \mathrm{E}(\Pr(F[-2a, 2a]^c > \varepsilon | \theta_1, \ldots, \theta_n) | X_1, \ldots, X_n).$$

It is well known (cf. [4]) that the conditional law of $F$ given $\theta_1, \ldots, \theta_n$ is the Dirichlet distribution with base measure $\alpha + \sum_{i=1}^n \delta_{\theta_i}$. In particular,

$$F[-2a, 2a]^c | \theta_1, \ldots, \theta_n$$

$$\sim \mathrm{beta}(\alpha[-2a, 2a]^c + N[-2a, 2a]^c, \alpha[-2a, 2a] + N[-2a, 2a]),$$

where $N(A) = \sum_{i=1}^n 1\{\theta_i \in A\}$. We can use the preceding display and Markov's inequality on the inner expectation on the right-hand side to see that

$$\Pr(F[-2a, 2a]^c > \varepsilon | X_1, \ldots, X_n)$$

$$\leq \frac{\alpha[-2a, 2a]^c + \sum_{i=1}^n \Pr(\theta_i \in [-2a, 2a]^c, \sigma \leq b\sigma_n | X_1, \ldots, X_n)}{\varepsilon(\alpha(\mathbb{R}) + n)}$$

$$+ \Pr(\sigma > b\sigma_n | X_1, \ldots, X_n).$$

Let $\theta_{-n} = (\theta_1, \ldots, \theta_{n-1})$, $H(\theta_1, \ldots, \theta_n)$ be the joint distribution of $(\theta_1, \ldots, \theta_n)$, $H_n(\theta_n|\theta_{-n})$ be the conditional distribution of $\theta_n$ given $\theta_{-n}$ and $H_{-n}(\theta_{-n})$ the marginal distribution of $\theta_{-n}$. Bayes' formula then gives

$$\Pr(\theta_n \in [-2a, 2a]^c, \sigma \leq b\sigma_n | X_1, \ldots, X_n)$$

$$= \frac{\int_0^b \int \int_{t_n \in [-2a, 2a]^c} \prod_{i=1}^n s^{-1} e^{-(X_i - t_i)^2/2s^2} \, dH_n(t_n|t_{-n}) \, dH_{-n}(t_{-n}) \, dG_n(s)}{\int \int \int \prod_{i=1}^n s^{-1} e^{-(X_i - t_i)^2/2s^2} \, dH_n(t_n|t_{-n}) \, dH_{-n}(t_{-n}) \, dG_n(s)}.$$

The conditional distribution $H_n(\cdot|\theta_{-n})$ of $\theta_n$ given $\theta_{-n}$ is

$$\theta_n|\theta_{-n} = \begin{cases} \theta_i, & \text{with probability } 1/(\alpha(\mathbb{R}) + n - 1), \, i = 1, \ldots, n - 1, \\ \sim \bar{\alpha}, & \text{with probability } \alpha(\mathbb{R})/(\alpha(\mathbb{R}) + n - 1). \end{cases}$$

Thus, with $\delta$ a lower bound for the density of $\bar{\alpha}$ on $[-a, a]$ and $s < a$,

$$\int e^{-(X_n - t_n)^2/2s^2} \, dH_n(t_n|t_{-n}) \geq \frac{\alpha(\mathbb{R})}{\alpha(\mathbb{R}) + n - 1} \int_{X_n - s}^{X_n + s} e^{-(X_n - t_n)^2/2s^2} \, d\bar{\alpha}(t_n)$$

$$\geq \frac{\alpha(\mathbb{R})}{\alpha(\mathbb{R}) + n - 1} e^{-1/2} \delta s,$$

provided that $s < a$. Thus, the integral in the denominator of the Bayes formula is bounded below by

$$(8.2) \quad \frac{\alpha(\mathbb{R})}{\alpha(\mathbb{R}) + n - 1} e^{-1/2} \delta \int_0^a \int \prod_{i=1}^{n-1} s^{-1} e^{-(X_i - t_i)^2/2s^2} \, dH_{-n}(t_{-n}) \, dG_n(s).$$

We now upper-bound the numerator. For $|X_n| \leq a$ and $t_n \in [-2a, 2a]^c$, we have that $(X_n - t_n)^2 \geq a^2$, so for any $s \leq b\sigma_n$ it follows that

$$s^{-1} e^{-(X_n - t_n)^2/2s^2} \leq s^{-1} e^{-a^2/4s^2} e^{-a^2/4b^2\sigma_n^2} \leq A_0 e^{-a^2/4b^2\sigma_n^2},$$

where $A_0 = \sup\{s^{-1} e^{-a^2/4s^2} : s > 0\}$. This leads to the bound

$$(8.3) \quad A_0 e^{-a^2/4b^2\sigma_n^2} \int_0^{b\sigma_n} \int \prod_{i=1}^{n-1} s^{-1} e^{-(X_i - t_i)^2/2s^2} \, dH_{-n}(t_{-n}) \, dG_n(s).$$

As $b\sigma_n < a$, the ratio of the integral in (8.3) over that in (8.2) is bounded by 1. Thus, we may bound the expression in the Bayes formula by $Kne^{-a^2/4b^2\sigma_n^2}$, for some constant $K$. Putting this into the bound for $\Pr(F[-2a, 2a]^c > \varepsilon | X_1, \ldots, X_n)$, we complete the proof of the first assertion.

For the second assertion, observe that the restriction $X_i \in [-a, a]$ is redundant for sufficiently large $a$. Replace $\varepsilon$ by $\varepsilon_n^2$, where $\varepsilon_n$ satisfies (4.2), and $b$ by a sequence $b_n$ that satisfies the given conditions. It then follows from Fubini's theorem (cf. the proof of Theorem 2.1 of [8]) that $\Pi(\sigma > b_n\sigma_n | X_1, \ldots, X_n) \to 0$ in $P_0^n$-probability. Since $n\varepsilon_n^2 \to \infty$, the result follows. $\square$

**9. Proof of Theorem 1.** We apply Theorem 5, with $\varepsilon_n$ given by (2.1), and

$$\mathcal{P}_n = \{p_{F,\sigma} : F[-a,a]^c \leq \varepsilon_n^2\},$$

$$\mathcal{P}_{n,j} = \{p_{F,\sigma} : F[-a,a]^c \leq \varepsilon_n^2, 2^j\sigma_n \leq \sigma < 2^{j+1}\sigma_n\}, \qquad j = 0, \pm 1, \ldots,$$

where $[-a/2, a/2]$ contains the support of $p_0$.

Let $\underline{\varepsilon}_n = (n\sigma_n)^{-1/2}\log n \vee \sigma_n^2 \log n$, which is smaller than $\varepsilon_n$. In the second part of the proof, it will be seen that (4.2) holds when $\underline{\varepsilon}_n$ replaces $\varepsilon_n$. If we choose $b_n$ to be a sufficiently large multiple of $(n\underline{\varepsilon}_n^2)^{1/\gamma}$, then $\Pr(\sigma > b_n\sigma_n) \lesssim e^{-\beta b_n^\gamma} \leq e^{-cn\underline{\varepsilon}_n^2}$ for an arbitrarily large constant $c$ and $b_n\sigma_n \leq \sigma_n^{1-1/\gamma}(\log n)^{2/\gamma} \vee (n\sigma_n^{4+\gamma}\log n)^{1/\gamma}$, which goes to 0 as a power of $n$ up to a log factor, since $n^{-a_1} \leq \sigma_n \leq n^{-a_2}$ and $a_2 \geq (4+\gamma)^{-1}$. Thus, all conditions of the second part of Lemma 11 hold and hence $\Pi_n(\mathcal{P}_n^c|X_1,\ldots,X_n) \to 0$ in probability.

By Lemma A.3 of [9], the $L_1$-distance between $p_{F,\sigma}$ and $p_{F',\sigma}$ for $F'$ equal to $F$ restricted and renormalized to $[-a,a]$ is bounded above by $2F[-a,a]^c$. Therefore, for $\varepsilon \geq \varepsilon_n$, we have, with $\mathcal{P}_{n,b_1,b_2}$ as in Lemma 3,

$$\log N(3\varepsilon^2, \mathcal{P}_{nj}, \|\cdot\|_1)$$

$$\leq \log N(\varepsilon^2, \mathcal{P}_{a,2^j\sigma_n, 2^{j+1}\sigma_n}, \|\cdot\|_1)$$

$$\lesssim \log\left(\frac{2}{\varepsilon^2}+1\right) + \left(\frac{a}{2^j\sigma_n}+1\right)\left(\log\frac{1}{\varepsilon^2}\right)\left(\log\left(\frac{a}{2^j\sigma_n}+1\right)\frac{1}{\varepsilon^2}\right),$$

for $\varepsilon^2 < 1/2$ by Lemma 3. Hence

$$\log N(\sqrt{3}\varepsilon, \mathcal{P}_{nj}, h) \leq \log N(3\varepsilon^2, \mathcal{P}_{nj}, \|\cdot\|_1)$$

$$\lesssim \begin{cases} 2^{-j}\dfrac{a}{\sigma_n}\left(\log\dfrac{a2^{-j}}{\varepsilon^2\sigma_n}\right)^2, & 2^j\sigma_n < a, \\ \left(\log\dfrac{1}{\varepsilon^2}\right)^2, & 2^j\sigma_n > a. \end{cases}$$

It follows that $N(\varepsilon_n, \mathcal{P}_{nj}, h)$ is bounded by a multiple of $\exp(Cj^2 2^{-j}\sigma_n^{-1}\log^2 n)$ if $2^j\sigma_n < a$, and is bounded by a multiple of $\exp(C\log^2 n)$ otherwise, for a large constant $C$. By the assumption on the prior of $\sigma$, we have

$$\Pi_n(\mathcal{P}_{nj}) \leq G(2^{j+1}) - G(2^j) \lesssim \begin{cases} e^{-\beta 2^{-\gamma(j+1)}}, & j < 0, \\ e^{-\beta 2^{\gamma j}}, & j \geq 0. \end{cases}$$

Thus, to verify (4.1) for a multiple of $\varepsilon_n$, it suffices to show that for constants $C, D, E$,

$$\sum_{j\leq 0}\exp(Cj^2 2^{-j}\sigma_n^{-1}\log^2 n - D2^{-\gamma j}) \lesssim e^{En\varepsilon_n^2},$$

$$\sum_{0\leq j\leq \log(a/\sigma_n)}\exp(Cj^2 2^{-j}\sigma_n^{-1}\log^2 n - D2^{\gamma j}) \lesssim e^{En\varepsilon_n^2},$$

$$\sum_{j \geq \log(a/\sigma_n)} \exp(C \log^2 n - D 2^{\gamma j}) \lesssim e^{En\varepsilon_n^2}.$$

For the third sum, this is immediate. In the second sum, we can bound the factors $j^2 2^{-j}$ by a constant and the inequality is then immediate.

The first sum can be transformed into a sum for $j = 0, 1, \dots$ by the change of variable $j \mapsto -j$. The factor $j^2$ can be absorbed into $2^j$ at the expense of replacing 2 by a slightly bigger number $A = 2^\eta$, where $\eta > 1$ can be arbitrarily close to 1. Put $S(K) = \exp(KA^j - A^{\gamma' j})$, where $\gamma' = \gamma/\eta$. To study the growth rate of $\log S(K)$ as $K \to \infty$, observe that $KA^j - A^{\gamma' j}$ is maximized near $j_0 = (\gamma' - 1)^{-1} \log_A(K/\gamma')$, leading to function value at most a multiple of $K^{\gamma'/(\gamma'-1)}$, that is, $K^{(\gamma/(\gamma-1))+}$. Since the series decays faster than geometrically, the sum in the tail is bounded by a multiple of the maximum term. The first $j_0$ terms together contribute at best $j_0$ times the value of the maximum. Thus, it follows that $\log S(K) = O(K^{(\gamma/(\gamma-1))+} \log K)$ and, clearly, the logarithmic factor can be absorbed into the power. Hence, in view of the fact that $\log \varepsilon_n^{-1} = O(\log n)$, the requirement (4.1) becomes $(a/\sigma_n)^{(\gamma/(\gamma-1))+} (\log \frac{a}{\varepsilon_n^b \sigma_n})^{(2\gamma/(\gamma-1))+} \lesssim n\varepsilon_n^2$. Again, the logarithmic factor may be absorbed into the power. Thus, the condition is satisfied in view of (2.1).

Finally, we verify (4.2). Fix numbers $b' > b > 0$, to be chosen sufficiently large at the end of the proof. Because $P_0$ possesses compact support, by Lemma 2, there exists a discrete distribution $F_n = \sum_{j=1}^{N_n} p_j \delta_{z_j}$, supported on $N_n \lesssim \sigma_n^{-1} \log \varepsilon_n^{-b'}$ points in the $\sigma_n$-enlargement of the support of $P_0$, such that

$$\tag{9.1} \|p_{F_n, \sigma_n} - p_{P_0, \sigma_n}\|_1 \lesssim \varepsilon_n^{b'} (\log \varepsilon_n^{-b'})^{1/2} \lesssim \varepsilon_n^b,$$

for sufficiently large $n$. This will change by $O(\varepsilon_n^b)$ if we move the support points of $F_n$ by $2\varepsilon_n^b \sigma_n$, so we can assume that the support points are $\varepsilon_n^b \sigma_n$-separated. We can then find disjoint intervals $U_1, \dots, U_{N_n}$ with $z_j \in U_j$ and $\lambda(U_j) = \varepsilon_n^b \sigma_n$ for $j = 1, \dots, N_n$. We can modify this to a partition of an interval that contains the support of $P_0$ into $M_n \geq N_n$ intervals $U_1, \dots, U_{M_n}$, such that each interval $U_j$ has length $\varepsilon_n^b \sigma_n \leq \lambda(U_j) \leq 2\varepsilon_n^b \sigma_n$ for $j = 1, \dots, M_n$, and such that $M_n \lesssim \sigma_n^{-1} \log n$.

Let $\sum_{j=1}^{N_n} |F(U_j) - p_j| \leq \varepsilon_n^b$ and $F(U_j) \geq \varepsilon_n^{2b}$ for $j = 1, \dots, M_n$, and $|\sigma - \bar{b}\sigma_n| \leq \varepsilon_n^b \sigma_n$, for $\bar{b} = (b_1 + b_2)/2$. By Lemma 5 [with $U_0 = (\bigcup_{1 \leq j \leq N_n} U_j)^c$], we have $\|p_{F, \sigma_n} - p_{F_n, \sigma_n}\|_1 \lesssim \varepsilon_n^b$. Furthermore, $\|p_{F, \sigma} - p_{F, \sigma_n}\|_1 \leq \|\phi_\sigma - \phi_{\sigma_n}\|_1 \leq |\sigma - \sigma_n|/(\sigma \wedge \sigma_n) \leq \varepsilon_n^b$. By the triangle inequality and (9.1), we have $\|p_{P_0, \sigma_n} - p_{F, \sigma}\|_1 \lesssim \varepsilon_n^b$ and hence $h(p_0, p_{F, \sigma}) \lesssim \varepsilon_n^{b/2}$. Combining the preceding inequality with Lemma 4, we conclude that $h(p_0, p_{F, \sigma}) \lesssim \sigma_n^2 + \varepsilon_n^{b/2} \lesssim \sigma_n^2$ if $b$ is sufficiently large. Now, if $x \in \text{supp}(p_0)$, then

$$p_{F, \sigma}(x) \geq \int_{x-\sigma_n}^{x+\sigma_n} \phi_\sigma(x-z) \, dF(z) \gtrsim e^{-(\sigma_n/\sigma)^2/2} \frac{F[x - \sigma_n, x + \sigma_n]}{\sigma} \gtrsim \frac{\varepsilon_n^{2b}}{\sigma_n},$$

because the interval $[x - \sigma_n, x + \sigma_n]$ contains at least one of the intervals $U_1, \ldots, U_{M_n}$. Consequently, $P_0(p_0/p_{F,\sigma}) \lesssim \sigma_n/\varepsilon_n^{2b}$.

Results of the two preceding paragraphs imply, for sufficiently large $b$, that

$$\left\{ (F, \sigma): \sum_{j=1}^{N_n} |F(U_j) - p_j| \leq \varepsilon_n^b, F(U_j) \geq \varepsilon_n^{2b}, j = 1, \ldots, M_n, |\sigma - \bar{b}\sigma_n| \leq \varepsilon_n^b \sigma_n \right\}$$

$$\subset \left\{ (F, \sigma): h(p_0, p_{F,\sigma}) \lesssim \sigma_n^2, P_0\left(\frac{p_0}{p_{F,\sigma}}\right) \lesssim \frac{\sigma_n}{\varepsilon_n^{2b}} \right\}.$$

The densities $p_{F,\sigma}$ with $(F, \sigma)$ as in the last set are contained in $B(p_0, c_5\sigma_n^2 \log n) \subset B(p_0, c_4\varepsilon_n)$ for a sufficiently large constant $c_4$, in view of Lemma 7.

By construction, $\lambda(U_j) \gtrsim \varepsilon_n^b \sigma_n$ and hence $\alpha(U_j) \gtrsim \varepsilon_n^b \sigma_n$ for every $j = 1, \ldots, M_n$. Furthermore, for sufficiently large $b$, we have $\varepsilon_n^b M_n \leq 1$. By Lemma 10 (with $p_j = 0$ for $N_n < j \leq M_n$ and a different constant $b$, as in the present proof), we conclude that the prior mass of the set $B(p_0, c_5\varepsilon_n)$ is bounded below by a multiple of $\varepsilon_n^b \exp(-cM_n \log \varepsilon_n^{-1}) \geq \exp(-c'\sigma_n^{-1}(\log n)^2)$, proving the first statement.

For the proof of the last statement of the theorem, we follow the same steps, but we redefine $\varepsilon_n$ by (2.2). The verification of (4.1) needs no changes, but we adapt the verification of (4.2) as follows. Fix $b' > b > 0$. Because $P_0$ possesses compact support, by Lemma 2, there exists a discrete distribution $F_n = \sum_{j=1}^{N_n} p_j \delta_{z_j}$ supported on $N_n \lesssim \sigma_n^{-1} \log \varepsilon_n^{-b'}$ points in the $\sigma_n$-enlargement of the support of $P_0$ such that (9.1) holds for sufficiently large $n$. The proof of Lemma 2 shows that we can satisfy $F_n(I_j) = P_0(I_j)$ for every interval $I_j$ in a covering of the support of $P_0$ by $M_n \lesssim \sigma_n^{-1}$ intervals of length $\sigma_n$. We can assume that the support points are $\varepsilon_n^b \sigma_n$-separated. We can then find disjoint intervals $U_1, \ldots, U_{N_n}$ with $z_j \in U_j$ and $\lambda(U_j) = \varepsilon_n^b \sigma_n$ for $j = 1, \ldots, N_n$, and such that each $U_j$ is contained in some interval $I_k$.

Suppose that $F$ is a probability measure satisfying $\sum_{j=1}^{N_n} |F(U_j) - p_j| \leq \sigma_n^a \varepsilon_n^b$ and that $\sigma$ is a number with $|\sigma - \bar{b}\sigma_n| \leq \varepsilon_n^b \sigma_n$ for $\bar{b} = (b_1 + b_2)/2$. As before, this implies that $h(p_{P_0,\sigma_n}, p_{F,\sigma}) \lesssim \varepsilon_n^{b/2}$. Moreover, for every $x \in \mathrm{supp}(p_0)$, $p_{F,\sigma}(x) \gtrsim \sigma_n^{-1} F[x - \sigma_n, x + \sigma_n] \gtrsim \sigma_n^{-1} \min_j F(I_j)$, because the interval $[x - \sigma_n, x + \sigma_n]$ contains at least one of the intervals $I_1, \ldots, I_{M_n}$. By construction, $F_n(I_j) = P_0(I_j)$, which is bounded below by a multiple of $\sigma_n^a$, by assumption, for every $j$. Hence,

$$F(I_j) \geq \sum_{i:U_i \subset I_j} F(U_i) \geq \sum_{i:U_i \subset I_j} p_i - \sigma_n^a \varepsilon_n^b = F_n(I_j) - \sigma_n^a \varepsilon_n^b \gtrsim \sigma_n^a.$$

Consequently, $p_{F,\sigma}(x) \gtrsim \sigma_n^{a-1}$ if $x \in \mathrm{supp}(p_0)$ and so $\|p_0/p_{F,\sigma}\|_\infty \lesssim \sigma_n^{1-a}$.

By Lemma 4, we have that $h(p_0, p_{P_0,\sigma_n}) \lesssim \sigma_n^2$. Therefore, since $p_0/p_{P_0,\sigma_n}$ is bounded by Lemma 6, we conclude by Lemma 8 that $P_0(\log(p_0/p_{P_0,\sigma_n}))^k \lesssim$ $\sigma_n^4$, $k = 1, 2$. With the help of Lemma 9, with $p = p_0$, $q = p_{P_0, \sigma_n}$ and $r = p_{F, \sigma}$, we see that $P_0(\log(p_0/p_{F, \sigma}))^k \lesssim \sigma_n^4 + \varepsilon_n^b + \varepsilon_n^{b/2} \sigma_n^{1-a} \lesssim \sigma_n^4$, $k = 1, 2$, for sufficiently large $b$. Combining the results of the three preceding paragraphs, we see that, for sufficiently large $b$,

$$\left\{ (F, \sigma) : \sum_{j=1}^{N_n} |F(U_j) - p_j| \leq \sigma_n^a \varepsilon_n^b, |\sigma - \bar{b} \sigma_n| \leq \varepsilon_n^b \sigma_n \right\}$$

$$\subset \left\{ (F, \sigma) : P_0 \left( \log \frac{p_0}{p_{F, \sigma}} \right)^k \lesssim \sigma_n^4, k = 1, 2 \right\} \subset B(p_0, c_5 \sigma_n^2).$$

By construction, $\lambda(U_j) = \varepsilon_n^b \sigma_n$ and hence $\alpha(U_j) \gtrsim \varepsilon_n^b \sigma_n$ for every $j$. By Lemma 10, we conclude that the prior mass of the set $B(p_0, c_5 \sigma_n^2)$ is bounded below by a multiple of $\varepsilon_n^b \exp(-c N_n \log \varepsilon_n^{-b}) = \exp(-c' \sigma_n^{-1} (\log n)^2)$, proving (2.2).

The validity of the final remark is clear from the proof, as there are only finitely many terms when $G$ is compactly supported.

## 10. Proof of Theorem 2.

Fix a smooth function $w : \mathbb{R} \to [0, 1]$ with support $[-2, 2]$ that is identically 1 on $[-1, 1]$ and let $w_n(x) = w(x/k_n)$ for $k_n = (\log n / c')^{1/\gamma}$ for some $c' < c$. Define new observations $\bar{X}_1, \dots, \bar{X}_{\bar{n}}$ from the original observations $X_1, \dots, X_n$ by rejecting each $X_i$ independently with probability $w_n(X_i)$. Because $P_0[-k_n, k_n]^c = o(n^{-1})$ by the tail assumption on $P_0$, the probability that some $X_i$ is rejected is actually $o(1)$ and hence the posterior distributions based on the new and the original observations are the same with probability tending to one. In particular, they have the same posterior rate of convergence. The new observations are a random sample from the density $p_n$ that is proportional to $p_0 w_n$. Because $|\int p_0 w_n \, d\lambda - 1| \leq P_0[-k_n, k_n]^c = o(n^{-1})$, we have that $h_k(p_n, p_0) = o(n^{-1})$ for every $k$. Hence, it suffices to show that the posterior based on the new observations concentrates at rate $\varepsilon_n$ around $p_n$.

We shall establish this by means of an obvious triangular array version of Theorem 2.1 of [8] [with the only difference being that we treat $\Pi_n(\mathcal{P}_n^c) \to 0$ in $P_0^n$-probability directly instead of through their condition (2.3)]. We verify this for $\mathcal{P}_n = \{p_{F, \sigma} : F[2k_n, 2k_n]^c \leq 2\varepsilon_n\}$ and $\varepsilon_n = \max\{(n\sigma_n)^{-1/2}(\log n)^{1+1/(2\gamma)}, \sigma_n^2 \log n\}$. We choose $w$ such that $\int (w'/w)^4 w \, d\lambda < \infty$ and $\int (w''/w)^2 w \, d\lambda < \infty$. Then $\int (p_n'/p_n)^4 p_n \, d\lambda = O(1)$ and $\int (p_n''/p_n)^2 p_n \, d\lambda = O(1)$, and hence $h^2(p_n, p_n * \phi_{\sigma_n}) = O(\sigma_n^2)$, by Lemma 4.

The verification of the entropy bound can proceed as before, except that we obtain an additional logarithmic factor by the dependence of $a_n$ on $n$, as follows:

$$\log N(3\varepsilon_n^2, \mathcal{P}_n, \|\cdot\|_1) \lesssim \frac{a_n}{\sigma_n} \left( \log \frac{a_n}{\varepsilon_n^2 \sigma_n} \right)^2 \lesssim \frac{1}{\sigma_n} (\log n)^{2+1/(2\gamma)} \leq n \varepsilon_n^2.$$

For the verification that $\Pi_n(\mathcal{P}_n^c | X_1, \ldots, X_n) \to 0$, we use Lemma 11, where it suffices that

$$(10.1) \qquad e^{-a_n^2/(8\sigma_n^2 b_2^2)} \frac{1}{\varepsilon_n} \frac{n}{\sigma_n \min_{|t| \le k_n} \alpha'(t)} \to 0.$$

This is certainly the case under the tail condition on $\alpha'$.

We adapt the verification of prior concentration rate in Kullback–Leibler neighborhoods as follows. Fix $b' > b > 0$. Because $P_n$ possesses support $[-2k_n, 2k_n]$, by Lemma 2, there exists a discrete distribution $F_n = \sum_{j=1}^{N_n} p_j \delta_{z_j}$ supported on $N_n \lesssim k_n \sigma_n^{-1} \log \varepsilon_n^{-b'}$ points in the interval $[-2k_n - \sigma_n, 2k_n + \sigma_n]$ such that

$$(10.2) \qquad \| p_{F_n, \sigma_n} - p_{P_n, \sigma_n} \|_1 \lesssim \varepsilon_n^{b'} (\log \varepsilon_n^{-b'})^{1/2} \lesssim \varepsilon_n^b,$$

for sufficiently large $n$. Because this distance changes by $O(\varepsilon_n^b)$ if we move the support points of $F_n$ by $2\varepsilon_n^b \sigma_n$, we can assume that the support points are $\varepsilon_n^b \sigma_n$-separated. We can then find disjoint intervals $U_1, \ldots, U_{N_n}$ with $z_j \in U_j$ and $\lambda(U_j) = \varepsilon_n^b \sigma_n$ for $j = 1, \ldots, N_n$. We can modify this to a partition of the interval $[-2k_n, 2k_n]$ into $M_n \ge N_n$ intervals $U_1, \ldots, U_{M_n}$ such that each interval $U_j$ has length contained in $[\varepsilon_n^b \sigma_n, \sigma_n]$ and such that $M_n \lesssim k_n \sigma_n^{-1} \log \varepsilon_n^{-b}$.

Let $F$ satisfy $\sum_{j=1}^{N_n} |F(U_j) - p_j| \le \varepsilon_n^b$ and $F(U_j) \ge \varepsilon_n^{2b}$ for $j = 1, \ldots, M_n$, and suppose that $|\sigma - \bar{b}\sigma_n| \le \varepsilon_n^b \sigma_n$ for $\bar{b} = (b_1 + b_2)/2$. By Lemma 5, we have $\| p_{F, \sigma_n} - p_{F_n, \sigma_n} \|_1 \lesssim \varepsilon_n^b$ and $\| p_{F, \sigma} - p_{F, \sigma_n} \|_1 \le \varepsilon_n^b$. Applying the triangle inequality repeatedly and combining the preceding two inequalities with (10.2), we find that $\| p_{P_n, \sigma_n} - p_{F, \sigma} \|_1 \lesssim \varepsilon_n^b$ and so $h(p_{P_n, \sigma_n}, p_{F, \sigma}) \varepsilon_n^{b/2}$. Combining the preceding inequality with Lemma 4, we conclude that $h(p_n, p_{F, \sigma}) \lesssim \sigma_n^2 + \varepsilon_n^{b/2} \lesssim \sigma_n^2$ if $b$ is sufficiently large.

For every $x$ in the interval $[-2k_n, 2k_n]$, we have

$$p_{F, \sigma}(x) \ge \int_{x - \sigma_n}^{x + \sigma_n} \phi_\sigma(x - z) \, dF(z) \gtrsim \frac{F[x - \sigma_n, x + \sigma_n]}{\sigma} \gtrsim \frac{\varepsilon_n^{2b}}{\sigma_n},$$

because the interval $[x - \sigma_n, x + \sigma_n]$ contains at least one of the intervals $U_1, \ldots, U_{M_n}$. Consequently, $P_n(p_n/p_{F, \sigma}) \lesssim \sigma_n/\varepsilon_n^{2b}$.

Combining the above results, we see that, for sufficiently large $b$,

$$\left\{ (F, \sigma) : \sum_{j=1}^{N_n} |F(U_j) - p_j| \le \varepsilon_n^b, F(U_j) \ge \varepsilon_n^{2b}, j = 1, \ldots, M_n, |\sigma_n - \sigma| \le \varepsilon_n^b \sigma_n \right\}$$

$$\subset \left\{ (F, \sigma) : h(p_n, p_{F, \sigma}) \lesssim \sigma_n^2, P_n\left( \frac{p_n}{p_{F, \sigma}} \right) \lesssim \frac{\sigma_n}{\varepsilon_n^{2b}} \right\}.$$

The densities $p_{F, \sigma}$, with $(F, \sigma)$ as in the last set, are contained in $B(p_n, c_5 \sigma_n^2 \log n)$ for a sufficiently large constant $c_5$ in view of Lemma 7.

By construction, $\lambda(U_j) \gtrsim \varepsilon_n^b \sigma_n$ and every $U_j$ is contained in the interval $[-2k_n, 2k_n]$. By the lower bound assumption on $\alpha'$, we see that $\alpha(U_j) \gtrsim \min_{|t| \leq k_n} \alpha'(t) \varepsilon_n^b \sigma_n \gtrsim n^{-e}$ for every $j = 1, \ldots, M_n$ and some $e > 0$. Furthermore, we have $\varepsilon_n^b M_n \leq 1$ if we choose $b$ sufficiently large. By Lemma 10, we conclude that the prior mass of the set $B(p_n, c_5 \varepsilon_n)$ is bounded below by a multiple of $\exp(-cM_n \log \varepsilon_n^{-1}) = \exp(-c' \sigma_n^{-1} (\log n)^{2+1/\gamma})$.

**11. Proofs of Theorems 3 and 4.** Let $\mathcal{P}_n$ be $\mathcal{P}_{a_n, \sigma_n}$ in the notation of Theorem 6. Hence, its bracketing entropy integral is bounded by

$$\int_0^{\varepsilon_n} \sqrt{\frac{a_n}{\sigma_n} \left( \log \frac{a_n}{\varepsilon \sigma_n} \right)^2} \, d\varepsilon \lesssim \sqrt{\frac{a_n}{\sigma_n}} \varepsilon_n \log \frac{a_n}{\sigma_n \varepsilon_n}.$$

For $\varepsilon_n = (a_n/(n\sigma_n))^{1/2} \log n$, this is bounded above by a multiple of $\sqrt{n}\varepsilon_n^2$. Because $p_0$ has compact support, $\mathcal{P}_n$ contains $p_{p_0, \sigma_n} = p_0 * \phi_{\sigma_n}$, at least if $b_1 < 1 < b_2$, as we shall assume for simplicity. By Lemma 6, the quotient $p_0/p_{p_0, \sigma_n}$ is bounded above uniformly in $n$. The distance $h(p_0, p_{p_0, \sigma_n})$ is of the order $O(\sigma_n^2)$, by assumption. Hence, Theorem 3 follows by an application of Theorem 4 in [16], or Theorem 3.4.4 in [14].

The sieve $\mathcal{P}_n$ given by (3.2) is equal to the set $\mathcal{P}_{A, \sigma_n}$ considered in Theorem 7. For the given function $A$, the function $\bar{A}$ in this theorem can be taken to be equal to a multiple of $1 - \Phi(ra)$ if $\delta \leq 1/2$ for some $0 < r < 1$, and equal to $A(ra)$ for some $r < 1$ if $\delta \geq 1/2$. Therefore, $\bar{A}^{-1}(\varepsilon) \lesssim (\log \varepsilon^{-1})^{(1 \vee 2\delta)/2}$ and, in view of Theorem 7, the bracketing integral of $\mathcal{P}_n$ is bounded by

$$\int_0^{\varepsilon_n} \sqrt{\frac{\bar{A}^{-1}(\sigma_n \varepsilon^2)}{\sigma_n} \left( \log \frac{\bar{A}^{-1}(\sigma_n \varepsilon^2)}{\varepsilon^2 \sigma_n} \right)^2} \, d\varepsilon \lesssim \sqrt{\frac{1}{\sigma_n}} \varepsilon_n \left( \log \frac{1}{\sigma_n \varepsilon_n} \right)^{1 + (1 \vee 2\delta)/4},$$

which is $O(\sqrt{n}\varepsilon_n^2)$ for $\varepsilon_n = (n\sigma_n)^{-1/2} (\log n)^{1+(1 \vee 2\delta)/4}$. The remainder of the proof can be completed as before.

Department of Statistics                     Department of Mathematics
North Carolina State University              Vrije Universiteit
2501 Founders Drive                          De Boelelaan 1081a
Raleigh, North Carolina 27695-8203           1081 HV Amsterdam
USA                                          The Netherlands
E-mail: ghosal@stat.ncsu.edu                 E-mail: aad@cs.vu.nl